\documentclass{amsart}

\usepackage{graphicx} 
\usepackage[foot]{amsaddr}
\usepackage{enumerate}
\usepackage{vem_ns}
\usepackage{xcolor}
\usepackage[hidelinks]{hyperref}
\usepackage{float}
\usepackage{mathtools}
\usepackage{comment}
\usepackage{subcaption}

\newcommand{\ms}[1]{{ \color{magenta}{[MS: #1]}}}

\title[VEM for the Smagorinsky model]{General Order Virtual Element Approximation for the Smagorinsky turbulence model}

\author{Stefano Berrone$^{1, *}$, Karol L. Cascavita$^{1, *}$, Enrique Delgado Ávila$^{2}$, Samuele Rubino$^{3}$, Maria Strazzullo$^{1, *}$ and Fabio Vicini$^{1, *}$}

\address{$^1$ Politecnico di Torino, Department of Mathematical Sciences ``Giuseppe Luigi Lagrange'', Corso Duca degli Abruzzi, 24, 10129, Turin, Italy.}
\address{$^*$ INdAM-GNCS group member.}
\address{$^2$ Departamento Matemática Aplicada II \& IMUS, Universidad de Sevilla, Avda. de los Descubrimientos s/n, 41092 Sevilla, Spain.}
\address{$^3$ Departamento EDAN \& IMUS, Universidad de Sevilla, Avda. Reina Mercedes s/n, 41012 Sevilla, Spain.}

\begin{document}

\begin{abstract}
In this paper, we investigate a Smagorinsky model in a virtual element framework to simulate convection-dominated Navier-Stokes equations. We conduct a two-dimensional numerical investigation to assess the performance of the general order virtual element approximation in this context. First, we examine numerically the convergence of the method with respect to the meshsize to certify the novel virtual element numerical discretization, which includes, for the first time, a discretization of the Smagorinsky term. Moreover, we present a numerical study of a lid-driven cavity for different Reynolds numbers (up to 10000) and meshes (uniform, anisotropic, and isotropic with hanging nodes). The results highlight the main advantage of using the virtual elements method in this context: the isotropic refinement with hanging nodes enhances the accuracy of the solution compared to the anisotropic mesh, uses fewer degrees of freedom with respect to the uniform mesh, and yields the most stable behavior in terms of convergence of the Newton solver.
\end{abstract}
\maketitle
\section{Introduction}
The solution of the Navier-Stokes equations (NSE) is an ubiquitous goal across a wide range of industrial and scientific fields. The main challenges arise in the convection-dominated regimes, where the meshsize should comply with the Kolmogorov scale to accurately capture the complex flow features \cite{kolmogorov1941local, kolmogorov1941dissipation}. However, this requirement yields prohibitively computational costs, and using coarser meshes is not a practical alternative, as it might result in inaccurate simulations with spurious numerical oscillations or failure of the Newton solver to converge. To address these problems, numerical stabilization has been classically employed: the NSE model is modified according to the chosen stabilization strategy, allowing the use of a coarser mesh while still achieving accurate results and yet reducing computational efforts. For a survey on numerical stabilization approaches, we refer to the monograph \cite{roos2008robust}. Moreover, the reader may refer to the following list for several works on numerical stabilization of convection-dominated flows successfully applied in various contexts \cite{Ivagnes2025,Strazzullo2025,Zoccolan2024,Zoccolan2025237}.
This contribution focuses on the Smagorinsky model, which can be traced back to Joseph Smagorinsky \cite{Smagorinsky1963} and was conceived for the numerical simulation of turbulence in weather forecasting.  It is known to be the basic Large Eddy Simulation (LES)
model \cite{chacon2014mathematical,sagaut2001large} to simulate turbulence and has been applied in numerous settings due to its simplicity \cite{berselli2006mathematics, germano1991dynamic, iliescu2003numerical, Khani_Waite_2015, lilly1992proposed, pare1992existence, rebollo2014numerical, Rozema2022, Winckelmans2001}.

Specifically, we propose and analyze a discrete Smagorinsky model in a virtual element framework.
Virtual Element Methods (VEM) were introduced in 2013 with the seminal works \cite{BeiraoBrezziCangianiEtAl2013} and \cite{BeiraoBrezziMariniEtAl2014} to solve linear diffusion problems. It is nowadays considered and popularized as an extension of the Finite Element Methods (FEM) to more general polytopal meshes. The cornerstone of this method is the ability to build shape functions for general polygons, which is a key aspect to constructing a robust method that accommodates arbitrary shapes and orders. The use of the \emph{virtual} label arises due to the presence of non-polynomial shape functions, whose point-wise values and closed-form expressions are not explicitly known. Computations are possible through local projection operators and suitable degrees of freedom. A non-exhaustive list of VEM applications across various fields follows
\cite{Beirao2015a, Benedetto2016c, Borio2021, borio2024generalordervirtualelement, Dassi2021b, Dassi2024150, Kumar2024, Xu2024}. Moreover, a comprehensive documentation discussing the research developed over the last decade was recently published \cite{VEM_Review2023}. 

\medskip
In the realm of fluid dynamics, virtual element discretizations were first devised for the Stokes problem \cite{BeiraoLovadinaVacca2017, MANZINI2022176} and later extended to the Navier-Stokes equation in \cite{BeiraoLovadinaVacca2018}.
In the latter, the authors introduced the enhancement condition in the divergence-free virtual setting. Namely, a condition that allows $L^2$-projectors to be computable for higher orders. The key novelty of this VEM family lies in the divergence-free property, which is usually only weakly preserved in conventional discretization techniques. 
In this work, we take advantage of both this important feature and the versatility of VE meshes, which have great potential to accurately solve fluid dynamics problems in convection-dominated settings.

We propose a standard convergence analysis of the VEM approximation to the exact solution of the Smagorinsky model to certify the novel VEM discretization.
Moreover, we test the remarkable properties of VEM in a lid-cavity problem with Reynolds numbers up to 10000, comparing three types of meshes: (i) uniform, (ii) anisotropic, and (iii) isotropic with hanging nodes, i.e., points between two consecutive aligned edges of the tessellation. 
The ease in using hanging nodes on general meshes is a peculiarity of polygonal approximations, and in our tests, the isotropic mesh with hanging nodes proved to be the best option in terms of accuracy, computational efficiency, and robustness with respect to the Newton solver convergence.
All numerical tests have been performed for two-dimensional spatial domains.

The main novelties of this contribution are as follows:
\begin{itemize}
\item the discretization of the Smagorinsky model via general order VEM approximation with numerical validation through several convergence tests;
\item {the exploitation of the advantageous general polygonal VEM meshes in convection-dominated regimes combined with Smagorinsky stabilization in comparison with other mesh choices representing the state of the art in the literature.} 
\end{itemize}
The rest of the paper is organized as follows. The model problem is presented in its continuous form in Section \ref{sec:problem}. Section \ref{sec:vem} is devoted to the virtual element setting. Finally, in Section \ref{sec:results}, we present the numerical results. We close the paper by drawing some conclusions in Section \ref{sec:conc}.

\section{Problem Formulation} \label{sec:problem}

We introduce the notation that will be used throughout the manuscript.
We use the notation $(\cdot, \cdot)_\mathcal{O}$ to denote the $L^2$-scalar product for any domain\footnote{Specifically, considering a discrete spatial domain $\Omega$, $\mathcal{O}$ may represent the whole domain, an element of the discretization, or even the boundary of the element.} $\mathcal{O} \subset \mathbb{R}^2$. Although the presented VEM methodology applies to a general framework that also accommodates $\mathcal{O} \subset \mathbb{R}^3$, we restrict the presentation to the bidimensional case. In addition, let  $\left | \; \cdot \; \right |_{\ell^2}$ represents the Euclidean norm for vectors, and, with a slight abuse of notation, let also denote the Frobenius norm defined for all $\boldsymbol{\gamma} \in \mathbb{R}^{d \times d}$ as $|\boldsymbol{\gamma}|_{\ell^2}\coloneqq(\boldsymbol{\gamma}: \boldsymbol{\gamma})^{\frac{1}{2}}$. We indicate with $(1:n)$ the ordered set of natural numbers from $1$ to any integer $n>0$. 
\del{Moreover, the tilde over a variable will denote a generic virtual function (e.g., $\tilde{q}$) or a generic virtual space (e.g., $\tilde{Q}$).}
Regarding physical notations, we use the dimensionless
Reynolds number $Re \coloneqq {\overline{U}{L}}/{\nu}$, where $\nu>0$, $\overline{U}$, and $L$ denote the kinematic viscosity, the characteristic velocity, and the characteristic length of the problem at hand, respectively. This dimensionless number describes the flow regime.  For example, a large Reynolds number indicates a convection-dominated flow and turbulent behavior.

\subsection{The incompressible Navier-Stokes Equations}
\label{sec:NSE}

We aim to find the velocity $\bu$ and the ratio between the fluid pressure and its density $p$ on the domain $\Omega \subset \mathbb R^{2}$, such that  
\begin{equation}
    \left\lbrace 
    \begin{alignedat}{2}
    - \nabla \cdot \left (  \nu \nabla \bu \right ) +
    {(\nabla \bu) \bu} + \nabla p  &= \bforce &\qquad& \text{in $\Omega$}, \\
    \nabla{\cdot} \bu &= 0 &\qquad& \text{in $\Omega$},\\
    \bu &= \bm{0}  &\qquad& \text{on $\Gamma_D$}, \\
    \normal \cdot [p \bI - \nu \GRAD \bu] &= \bm{0}  &\qquad& \text{on $\Gamma_N $},
    \end{alignedat} \label{eq:NSE:strongform}
    \right. 
\end{equation}

with $\mathbf f$ an external forcing term per unit mass \cite{quarteroni2008numerical}. For the sake of simplicity, we refer to $p$ and $\mathbf f$ as the pressure variable and the forcing term, respectively.
We apply homogeneous Dirichlet and Neumann boundary conditions in the distinct portions of the domain boundary $\partial \Omega$, denoted $\Gamma_D$ and $\Gamma_N$, respectively, with $\Gamma_D \neq \emptyset$ and ${\Gamma_N} \coloneqq \partial \Omega\setminus \Gamma_N$. For simplicity, we will consider $\Gamma_D = \partial \Omega$ and $\Gamma_N = \emptyset$ from now on. Naturally, other possible boundary conditions can be employed.
We define the velocity and pressure spaces as 
\begin{equation*}
\bbU \coloneqq [H^1_{0}(\Omega)]^2 \quad \text{ and } \quad \mathbb Q \coloneqq L^2_0(\Omega) = \{q \in L^2(\Omega) \; | \, (q,1)_\Omega = 0\},  
\end{equation*}
and equipped with the norm $\|\boldsymbol v \|_{\bbU} \coloneqq |\boldsymbol v|_{[H^1(\Omega)]^2}$ for all $ \boldsymbol v\in \bbU$ and $\|q \|_{\mathbb Q} \coloneqq \|q\|_{L^2(\Omega)}$ for all $q\in \mathbb Q$. We remark that if $\Gamma_N \neq \emptyset$, then $\mathbb Q \coloneqq L^2(\Omega)$.
Thus, the weak formulation of the model problem \eqref{eq:NSE:strongform} 
reads as follows: find $(\bu, p) \in \mathbb U \times \mathbb Q$ such that:

\begin{equation}
        \left\lbrace 
        \begin{alignedat}{2}
        \nu a(\bu, \bv) + c(\bu,\bu,\bv) + b( \bv, p)  =& (\bforce,\bv)_\Omega &\qquad&  \forall \bv \in \mathbb U , \\
        b(\bu,q) =& 0 &\qquad& \forall q \in \mathbb Q,
        \end{alignedat} \label{eq:navier_stokes:weakform}
        \right. 
\end{equation}
with $\bforce \in [L^2(\Omega)]^2$, and with bilinears form $a: \mathbb U \times \mathbb U  \rightarrow \Real$, $b: \mathbb U   \times \mathbb Q \rightarrow \Real$, and trilinear form $c:  \mathbb U\times \mathbb U  \times \mathbb U \rightarrow \Real$   defined as 
\begin{eqnarray*}
    a(\bw, \bv) &\coloneqq& (\nabla \bw, \nabla \bv)_\Omega, \\
    b(\bw,q) &\coloneqq& -(q, \nabla \cdot \chg{\bv}{\bw})_\Omega, \\
    c(\bw,\bz,\bv) &\coloneqq& ((\GRAD \bz)\bw, \bv)_\Omega,
\end{eqnarray*}   
for every $\bw, \bv, \bz \in \mathbb U$ and $q\in \mathbb Q$. 

\subsection{The Smagorinsky turbulence model}
\label{sec:Smago}
The Smagorinsky model introduces an artificial diffusion to the NSE, acting as a general stabilizing term for convection-dominated flows.
The key feature of this model is that it approximates the stresses of the unresolved subgrid scales (due to turbulence) by a relation between the strain-rate and an Eddy viscosity.  This model is inherently discrete, and we delay the discussion of the discrete Smagorinsky model until Section \ref{sec:Smago:discrete}. Building up on the ideas in \cite{EnriqueROM}, the Smagorinsky model can be reinterpreted as the discretized version of the following continuous problem: given a spatial domain $\Omega$, find $(\bu, p)$ such that
\begin{equation}
    \left\lbrace 
    \begin{alignedat}{2}
    - \nabla \cdot \left (  \left (  \nu +  \nuSmago{\bu} \right ) \nabla \bu \right ) + {(\nabla \bu) \bu} + \nabla p  &= \bm{f} &\qquad& \text{in $\Omega$}, \\
    \nabla{\cdot} \bu &= 0 &\qquad& \text{in $\Omega$},\\
    \bu &= \bm{0}  &\qquad& \text{on $\partial \Omega$}. \\
    \end{alignedat} \label{eq:Smago:strongform}
    \right. 
\end{equation}
Namely, compared to the strong formulation of NSE in \eqref{eq:NSE:strongform}, an Eddy viscosity term $\nuSmago{\bu}$ is added. Let $\mathcal{T}$ denote a tessellation of the domain $\Omega$, decomposed in elements denoted by  $\cell$ with diameter $\hCell$.
Thus, the Eddy viscosity term in the Smagorinsky model is defined as $\nu_S : \U \rightarrow \Real$ such that
\begin{equation}
    \label{eq:Eddyviscosity}
    \nuSmago{\bu} = C_S^2 \sum_{\cell \in \mesh} \hCell^2 \left | \nabla \bu _{|_\cell}\right |_{\ell^2} \chi_{\cell},
\end{equation}
where $\chi_\cell$ is the characteristic function over the tessellation element, while $C_S$ is the Smagorinsky constant \cite{SamueleThesis}, which is classically set (guided by numerical experiments \cite{lilly1966representation}) as $C_S = 0.1$: in this paper we follow this standard choice. The Smagorinsky term \eqref{eq:Eddyviscosity} is strictly related to the discretization of the space. We remark that for $h_T \rightarrow 0$, the Eddy viscosity term $\nu_S(\bu) \rightarrow 0$. This totally complies with the Kolmogorov scale turbulence description: a refined mesh can capture the complex behavior of the flow and does not need stabilization, which is needed for a larger meshsize.
   In the next section, we specify the approximation by employing the VEM numerical discretization.
\section{VEM approximation}\label{sec:vem}

We briefly present the enhanced formulation for the divergence-free VEM introduced in \cite{BeiraoLovadinaVacca2018}. We first give some preliminaries regarding the mesh assumptions and the discrete setting.

\subsection{Mesh assumptions}
Let $({\mesh})_h$ denote a $h$-refined mesh sequence of $\Omega$, with $\mesh$ being a decomposition into a tessellation with mesh elements denoted by $\cell$. The boundary and faces of the element are denoted by $\dCell$ and $F$, respectively. {We remark that in this paper, we use the word face to denote the edges of the polygon.}   
The notation $\hCell$ and $\hFace$ stands for the elemental and the face diameter, respectively.  We define the global meshsize $h\coloneqq\max \limits_{\cell \in \mesh}{\hCell}$. 
Furthermore, we assume that there exists a real number $\rho >0$, such that:
\begin{enumerate}
    \item for all $\mesh$ in the sequence, for all $\cell \in \mesh$, $\cell$ is star-shaped with respect to a ball of radius $r_B \geq \rho \hCell$ (star-shaped property);
    \item for all $\cell \in \mesh$, and for any face $\face$ of $\cell$, $\hFace$ is comparable with $\hCell$ as $ \rho \hCell \leq \hFace \leq \hCell$ (shape regularity property).
\end{enumerate}
Alternatively, relaxed mesh assumptions can be adopted; we refer the reader to~\cite{Brenner2018, Sorgente2021} for further details.

\subsection{Polynomial decomposition} Let $\Pol^k(\mathcal{O})$ denote the space of polynomials on a domain $\mathcal{O}$ of degree up to $k$. 
A key idea in devising projectors in VEM for the vector formulation is to decompose vector polynomials over a cell $T$ into a gradient component and their orthogonal component. Thus, polynomial decomposition reads as
\begin{equation*}
    {[\PolT]^2 = \PolDecompGrad[k]   \oplus \PolDecompPerp[k]}, 
\end{equation*}
where  ${{\PolDecompGrad[k] \coloneqq \nabla \PolT[k+1]}}$ and  ${{\PolDecompPerp[k] \coloneqq \bx^\perp[\PolT[k-1]]}}$ with  $\bx^\perp \coloneqq (y, -x)$ \cite{DassiVacca2020}. Namely, by the multiplication action of $\bx^\perp$, we build the orthogonal component to the gradients. 
The scalar scaled monomial $m_{\bm \alpha} \in \PolT$ with $|\bm \alpha|_{\ell^2} \leq k$ is defined as: 
\begin{equation}
    m_{\bm \alpha} \coloneqq \prod_{i=1}^2 \left(  \frac{x_i - x^T_i}{\hCell} \right)^{\alpha_i},
    \label{eq:monom:scalar}
\end{equation}
where we use the classic notation of the multi-index $\bm \alpha \coloneqq (\alpha_1, \alpha_2)$ and we denote with $x_i^{\cell}$ the $i$-th coordinate of the barycenter of $\cell$. 
Vector scaled monomials are defined as $\vmonom_{(\bm \alpha_x, \bm \alpha_y)} \coloneqq (m_{\bm \alpha_x}, m_{\bm \alpha_y}) \in \vPolT$. 

\begin{remark}[Vector monomial decomposition and scaling]
    \emergencystretch=3em
    Explicit closed forms of the polynomial decomposition can be obtained for vector scale monomials, see \cite{DassiVacca2020} for a complete presentation.   For the sake of completeness, we recall the 2D case: 
    let $\bm \calB \coloneqq (\beta_x, \beta_y)$, be a vector gathering the coefficients $\beta_x, \beta_y$ for each direction, then
    \begin{align*}
        \vmonom_{(\bm \calB, \bm 0)} &= \frac{\hCell}{(|\bm \calB|_{\ell^2} +1)} \nabla m_{(\beta_x +1, \beta_y)} + \frac{\beta_y}{(|\bm \calB|_{\ell^2} +1)} \bx^\perp m_{(\beta_x, \beta_y -1)}, \\
        \vmonom_{(\bm 0, \bm \calB)} &= \frac{\hCell}{(|\bm \calB|_{\ell^2} +1)} \nabla m_{(\beta_x, \beta_y+1)} - \frac{\beta_x}{(|\bm \calB|_{\ell^2} +1)} \bx^\perp m_{(\beta_x-1, \beta_y)}.
        \label{eq:monom:vector:decomposition}
    \end{align*}
    See also \cite{DassiVacca2020} for detailed guidelines on how to build VEM.
    Moreover, a recent review on the proper scaling of normalized monomials and its influence on the conditioning was reported in \cite{Cicuttin2025}.  
\end{remark}
\emergencystretch=3em
A monomial scalar and vector basis are denoted as $\Mspace\coloneqq\left\lbrace  \smonom_{\bm \alpha} : |\bm \alpha|_{\ell^2} \leq k  \right\rbrace$  and $\vMspace$, respectively.

\subsection{Local spaces}

Let a polynomial of degree $k\geq 2$ be fixed. The local virtual element space introduced originally in \cite{BeiraoLovadinaVacca2017} is \\
\begin{equation*}
    \begin{alignedat}{2}
        \Uspace^k_\cell \coloneqq \{ \bv \in [H^1(\cell)]^2 \; \vert           (i)& \quad{\bv_{|\dCell}}\
                 {\in [\Bspace]^2},
                 \\                                  
                 (ii)&\quad {{- \Delta \bv - \nabla q}}
                {{\ \in \PolDecompPerp[k]}}  \text{ for some } q \in L^2(T),\\
                 (iii)& \quad\Div{\bv}\
                 {{\in \Pol^{k-1}(\cell)}} 
        \},
    \end{alignedat}
\end{equation*}
where the first condition in the virtual space concerns the conformity prescription on the boundary with
\begin{equation*}
     {{\Bspace \coloneqq \{v \in C^0(\partial \cell) | v_{|\face} \in \Pol^{k}(\face), \forall \face \subset\dCell\}}}.
\end{equation*}
Then, the enhanced local virtual space introduced in \cite{BeiraoLovadinaVacca2018} is defined as
\begin{equation}
    \Vcell \coloneqq \{ \chg{\vhv}{\bv} \in \Uspace^k_\cell | (\chg{\vhv}{\bv} - \ProjNabla[k]{\chg{\vhv}{\bv}}, \bg^\perp)_\cell = 0, \forall \bg^\perp \in \PolDecompPerp[k]\setminus \PolDecompPerp[k-2] \}.    
\end{equation}
The enhancement uses an additional condition, \add{formulated in terms of the VE projector $\ProjNablaOper[k]$}, that enables the computation of higher-order projectors, e.g., the $L^2$-projector of order $k$, \add{$\ProjZero[k]{\vhv}$}. \add{The  definition of the local projectors is postponed to Section \ref{sec:proj}}. 

Let $\{\vbasisv^i\}_{1\leq i \leq dim(\Vcell) }$ be a basis of $\Vcell$ made of virtual functions which are not explicitly known. One property of $\Vcell$ is the polynomial inclusion $\vPolT \subseteq \Vcell$. 
\subsection{Degrees of freedom (DoFs)}
Let $\Ndofs\coloneqq dim(\chg{[\PolT]^2}{\Vcell})$ denote the local number of degrees of freedom. Let $N^{\add{v}}$ denote the number of vertices in a polygon $\cell$ \add{and $\NVertex=2N^v$ be the number of DoFs related to the vertices.} Let $\NdofEdge \coloneqq \add{2N^v} (k-1)$ denote the number of DoFs set on $\dCell$, excluding vertices, let $\NdofBoundary \coloneqq \NVertex+\NdofEdge$ denote the total DoFs on $\dCell$. Moreover, let $\NdofDiv$ and  $\NdofLap$ denote the number of DoFs related to the divergence and the Laplacian operators, respectively. The latter two are integrals evaluated inside the element acting as DoFs.  The velocity local degrees of freedom function $\dofOper[\cell]{\mathlcal{l}}:  \Vcell \rightarrow \Real^{\Ndofs}$ is defined as: 
\begin{equation*}
    \dof[\cell]{\mathlcal{l}}{\vhv} \coloneqq 
        \left\lbrace 
    \begin{alignedat}{4}
        &\del{\dofV{i,\add{\mathlcal{d}}}{\vhv} \coloneqq} \vhv(\bx^\Vertex_{\cell,i})_{\add{\mathlcal{d}}}, & \qquad & i \in (1 : \add{\tfrac{1}{2}}\NVertex ) \; ,\mathlcal{l} = \add{\tfrac{1}{2}\NVertex(\mathlcal{d}-1)+} i, \\   
        &\del{\dofE{i,\add{\mathlcal{d}}}{\vhv} \coloneqq} \vhv(\bx^\Edge_{\cell,i})_{\add{\mathlcal{d}}}, & &  i \in (1 :\add{\tfrac{1}{2}}\NdofEdge), \; \mathlcal{l} = \add{{\tfrac{1}{2}\NdofEdge(\mathlcal{d}-1)} + } \NVertex + i ,  \\ 
        &\del{\dofD{i}{\vhv}\coloneqq} (\Div{\vhv},\smonom_i)_\cell, && i \in (1 : \NdofDiv), \; \mathlcal{l}= \NdofBoundary + i, \\
        &\del{\dofL{i}{\vhv}\coloneqq}(\vhv, \vmonom_i)_\cell,  & & i \in (1:\NdofLap), \; \mathlcal{l} = \NdofBoundary  + \NdofDiv +i, 
    \end{alignedat}
    \right.
\end{equation*}
\add{ where $\add{\mathlcal{d}}\in\{1,2\}$  denotes the vector component index}, for all $\smonom_i \in \Mspace[k-1]\setminus \{1\}$ and $i \in (1 : \NdofDiv)$, for all $\vmonom_i \in  \bx^\perp\Mspace[k-3]$ and $i \in (1 : \NdofLap )$, and for all $\vhv \in \Vcell$ and $\mathlcal{l} \in (1 : \Ndofs)$, with $\bx^\Vertex_{\cell,i}$ the $i$-th vertex of the polygon $\cell$ and $\bx^\Edge_{\cell,i}$ the $i$-th edge point, usually a Gauss-Lobatto quadrature point (excluding vertices) located at the edge.  
The advantage of directly using the values at the quadrature rule points as DoFs is the ease of computing the integrals. An important feature of the virtual basis is that $\dof[\cell]{\mathlcal{l}}{\vbasisv^j} = \delta_{\mathlcal{l}j}$. 

\subsection{\add{Local} projectors}\label{sec:proj}
For all $\cell \in \mesh$, the $H^1$-projection operator $\ProjNablaOper[k]: \Vcell \rightarrow \vPolT$ is defined such that, for all $\vhv \in \Vcell$, 
\begin{equation} 
    \left\lbrace
    \begin{alignedat}{3}
        (\nabla \ProjNabla[k]{\vhv}, \nabla \bp)_\cell &= (\nabla \vhv, \nabla \bp)_\cell &
     \qquad \forall \ \bp \in \vPolT \setminus \vPolT[0],  \label{eq: proj_nabla} \\
        (\ProjNabla[k]{\vhv},  \bp)_\dCell &= (\vhv,  \bp)_\dCell & \qquad \forall \ \bp \in \vPolT[0].
    \end{alignedat} \right. 
\end{equation}
In a similar fashion, the $L^2$-projection operator $\ProjZero[k]~: \Vcell \rightarrow \vPolT$ is obtained by solving the following problem:  for all $\vhv \add{\in \Vcell}$,  
\begin{align}
    ( \ProjZero[k]{\vhv}, \bp)_\cell = (\vhv, \bp)_\cell & \qquad \forall \ \bp \in \vPolT.  
    \label{eq: proj_zero} 
\end{align} 
Finally, we define the $L^2$-projector for tensors, specifically for the velocity gradient, $\ProjZeroGradOper[k-1]:  \Vcell \rightarrow \PolTensorSpace$ such that, for all $\vhv \in \Vcell$ 
\begin{equation}
    (\ProjZeroGradOper[k-1]{\vhv}, \bm \tau)_\cell = (\nabla \vhv, \bm \tau)_\cell \qquad \forall \bm \tau \in \PolTensorSpace. 
\end{equation}

\subsection{Discrete global weak formulation} 

Let $k>1$ be the order of  VEM, then the global enhanced spaces are defined as follows
    \begin{equation*}
        \begin{alignedat}{2}
            \Vh &\coloneqq \{ \vhv \in [H^1_0(\Omega)]^2| {\vhvT} \in \Vcell, \forall \cell \in \mesh\}, \\
            \PressureSpaceglobal &\coloneqq \{ \qh \in L^2_0(\Omega)| q_{h|\cell} \in \PressureSpacelocalDef, \forall \cell \in \mesh\}.
        \end{alignedat}            
    \end{equation*}
Notice that the conformity prescription is only applied to the discrete velocity, while the discrete global pressure is a piecewise polynomial that can jump.  
The discrete version of the weak problem reads as: find $(\uhv, \ph) \in \Vh \times \PressureSpaceglobal $
such that
\begin{equation}
    \left\lbrace 
    \begin{alignedat}{2}
    \nu \Ah{\uhv}{\vhv} + {c_h(\uhv,\uhv,\vhv)}+ b_h(\vhv, \ph)  =&\   l_h(\vhv), &\qquad&  \forall \vhv \in \Vh , \\
    b_h(\uhv,\qh) =& \ 0, &\qquad& \forall \qh \in \PressureSpaceglobal,
    \end{alignedat} \label{eq:navier:weakform:discrete}
    \right. 
\end{equation}
for the global linear form $l_h: \Vh \rightarrow \Real$ formulated as
\[
        \displaystyle l_h(\vhv) \coloneqq \sum_{\cell \in \mesh}(\bforce,{\ProjZero[k]{\vhv}})_\cell,
\]
and the global bilinear forms $a_h: \Vh \times \Vh \rightarrow \Real$ and $b_h: \Vh \times \PressureSpaceglobal \rightarrow \Real$ consisting of the collection of local contributions as 
    \begin{align*}
        \displaystyle \Ah{\whv}{\vhv}  &\coloneqq \sum_{\cell \in \mesh} (\GRAD \ProjNabla[k]{\whv}, \GRAD \ProjNabla[k]{\vhv})_\cell + \Sh{\whv}{ \vhv},
        \\
        \displaystyle b_h(\vhv, \qh) &\coloneqq -\sum_{\cell \in \mesh} (\qh, \Div{\vhv})_\cell.
    \end{align*}    
The stabilization term $\ShOper$ is needed to ensure coercivity in the discrete system. We use the standard {dofi-dofi} stabilization $\ShOper: \Vh \times \Vh \rightarrow \Real$ defined as
        \[   \Sh{\whv}{\vhv} \coloneqq \sum_{\cell \in \mesh} \sum_{\ell}^{\Ndofs}  \dof[\cell]{\ell}{(\bI - \ProjNablaOper[k]) \whv}  \dof[\cell]{\ell}{(\bI - \ProjNablaOper[k]) \vhv}, \]
where $\bI$ denotes the identity matrix.  The trilinear form related with the nonlinear term $c_h:  \Vh\times \Vh \times \Vh \rightarrow \Real$ is defined such that, for all $(\whv, \xhv, \vhv) \in \Vh\times \Vh \times \Vh$,
\begin{equation}
    c_h(\whv\chg{;}{,}\xhv,\vhv) \coloneqq \sum_{\cell \in \mesh} (({\ProjZeroGrad[k-1]{\xhv}}) \ProjZero[k]{\whv}, \ProjZero[k]{\vhv})_\cell.    
\end{equation}
Different configurations of linear forms can arise from the various available projectors. We use here the original setting presented in \cite{BeiraoLovadinaVacca2018}.
\subsection{The discrete Smagorinsky model}\label{sec:Smago:discrete} The discretized Smagorinsky form $\AsmagoOper: \Vh \times \Vh \times \Vh \rightarrow \Real$, in a VEM framework, is written as follows
\begin{eqnarray*}
        \displaystyle \Ahsmago{\whv}{\xhv}{\vhv} &\coloneqq& \sum_{\cell \in \mesh} (\nuhSmago{\whv}\GRAD \ProjNabla[k]{\xhv}, \GRAD \ProjNabla[k]{\vhv})_\cell,   
\end{eqnarray*} 
where we have also introduced the discrete version of the Smagorinsky viscosity adapted to VEM and defined as $\nuhSmagoOper: \Vh \rightarrow \Real$ such that 
\begin{equation}
\label{eq:smagoVemterm}
    \nuhSmago{\whv} \coloneqq \Cs^2\sum_{\cell \in \mesh} \hCell^2 | \ProjZeroGrad{\whv}|_{\ell^2}\reviewerB{}{ \chi_{\cell}}.
\end{equation}
The global problem, including the Smagorinsky diffusive term, has the following form: 
find $(\uhv, \ph) \in \Vh \times \PressureSpaceglobal $
such that
\begin{equation}
    \left\lbrace 
    \begin{alignedat}{2}
    & \nu \Ah{\uhv}{\vhv} + \Ahsmago{\uhv}{\uhv}{\vhv}  \\
    & \qquad \qquad  \qquad \quad +  {c_h(\uhv,\uhv,\vhv)}+ b_h(\vhv, \ph)  =   l_h(\vhv), &\qquad&  \forall \vhv \in \Vh , \\
    &  b_h(\uhv,\qh) = \ 0, &\qquad& \forall \qh \in \PressureSpaceglobal.
    \end{alignedat} \label{eq:navier:weakform:discreteSmago}
    \right. 
\end{equation}
We use the shorthand notation $\bm{\mathcal{F}}_h^{NS}(\uhv, \ph)=\mathbf{0}$ for the Navier-Stokes problem \eqref{eq:navier:weakform:discrete} and $\bm{\mathcal{F}}_h(\uhv, \ph)=\mathbf{0}$ for the global problem including the Smagorinsky term \eqref{eq:navier:weakform:discreteSmago}. To solve the nonlinear global problem, we use the Newton scheme detailed in Algorithm~\ref{alg:discrete:newton_smago}. 
For completeness, we write the directional derivative of the Smagorinsky term along the direction $\delta \uhv$, readily verified to be of the following form:
\begin{eqnarray*}
    \partial_{\uhv} \displaystyle \Ahsmago{\uhv}{\uhv}{\vhv} [\delta \uhv]&=& 
     \sum_{\cell \in \mesh} \left(\nuhSmago{\uhv} \GRAD \ProjNabla[k]{\delta \uhv}, \GRAD \ProjNabla[k]{\vhv} \right)_\cell \\
    &&
     +\Cs^2 \sum_{\cell \in \mesh} \hCell^2 \left(t,(\GRAD \ProjNabla[k]{\uhv}: \GRAD \ProjNabla[k]{\vhv} )\right)_\cell,
\end{eqnarray*}
with $t = (\ProjZeroGrad{\uhv} : \ProjZeroGrad{\delta \uhv})/|\ProjZeroGrad{\uhv}|_{\ell^2}$.

\begin{algorithm}[h]
\caption{The Newton solver for the Smagorinsky model}\label{alg:discrete:newton_smago}
\begin{algorithmic}[1]
\State {Choose $\bU_h^0\coloneqq(\uhv^0, p_h^0)$}
\State {Choose $\epsilon >0$, set $\chg{(\bm R^0)_i= \infty}{\bm R_h= \bm \infty} $ and $\niter=0$} 
\While{$|| \chg{\bm R_h^\niter}{\bm R_h}|| > \epsilon $}
    \State {Seek $\bm \delta_h \coloneqq(\delta \uhv, \delta p_h)$ solving  for all  $(\vhv, q_h) \in \Vh \times \PressureSpaceglobal$ such that
    \[\partial_{\bU_h} \bm{\mathcal{F}}_h[\bm \delta_h] = -\bm{\mathcal{F}}_h, \]
    \hspace{0.5 cm} or explicitly written
\[\begin{cases}
    \partial_{\uhv} \bm{\mathcal{F}}_h^{NS}[\bm \delta_h]   + \partial_{\uhv} \Ahsmago{\uhv^\niter}{\uhv^\niter}{\vhv}[\bm \delta_h]  &= -R_1^s(\bU_h^\niter),\\
    \partial_{p_h} \bm{\mathcal{F}}_h^{NS}[\bm \delta_h] &=-R_2(\bU_h^\niter),
\end{cases}
\]
 \hspace{0.5 cm} where the residual is defined as $\bm R_h \coloneqq(R_1^s(\bU_h), R_2(\bU_h))$ and
\begin{eqnarray*}
    R_1^s(\bU_h) &\coloneqq&  R_1(\bU_h) + \chg{a_h^s}{\AsmagoOper}(\uhv; \uhv,\vhv),\\
     R_1(\bU_h) &\coloneqq&  c_h(\uhv, \uhv, \vhv ) +a_h(\uhv, \vhv) + b_h(\vhv,\ph) -  l_h(\vhv),   \\
     R_2(\bU_h) &\coloneqq&  b_h(\uhv,\qh).
\end{eqnarray*}
}
\State Update the solution $\bU_h^{\niter+1}= \bU_h^\niter + \bm \delta_h$ 
\State $\niter = \niter+1$
\EndWhile
\end{algorithmic}
\end{algorithm}

\section{Numerical Results}\label{sec:results}

In this section, we illustrate numerical results for the Smagorinsky model combined with the VEM approximation. We will show that the Smagorinsky model is necessary to simulate these Reynolds numbers on coarse meshes. Indeed, we conducted some numerical tests solving the standard NSE using either VEM (employing an in-house code based on \cite{berrone2025polydim}) or FEM (employing the FEniCS library \cite{fenics}) and \reviewerA{the simulation stalled}{the Newton solver did not converge} in the Reynolds number range $1000 \leq Re \leq 2000$. \\ 
We consider first some test cases assessing the accuracy of the model against manufactured or analytical solutions, and second, we consider the benchmark of the lid-driven cavity for increasing Reynolds numbers (up to $Re=10000$).
All problems are solved in the unitary square $\Omega$. The following metrics will be used for the computation of the errors for \reviewerAddA{the virtual element approximation of} velocity and pressure
\begin{align*}
\|\nabla \bu - \ProjZeroGrad[k-1]{\uhv}\|_{\mesh} & \coloneqq \left ( \sum_{\cell \in \mesh} \| \nabla \bu - \ProjZeroGrad[k-1]{\uhv} \|^2_{[L^2(\cell)]^{2\times 2}}\right )^{\frac{1}{2}}, \\
\|\bu - \ProjZero[k]{\uhv}\|_{\mesh} & \coloneqq \left ( \sum_{\cell \in \mesh} \| \bu - \ProjZero[k]{\uhv} \|^2_{[L^2(\cell)]^2}\right )^{\frac{1}{2}}, \\\
\| p - p_h \|_{\mesh} & \coloneqq \left( \sum_{\cell \in \mesh} \| p - p_h \|^2_{L^2(\cell)}
\right)^{\frac{1}{2}}.
\end{align*}

\subsection{\reviewerAddB{}{Accuracy} tests}
We consider two test cases here labeled as the irrotational force test and the polynomial P2P1 test. Let us detail the data of the two experiments: 
\begin{enumerate}
    \item \emph{Irrotational force test}: let $(\bu, p) \in [\Pol^2(\Omega)]^2 \times \Pol^{3}(\Omega)$ with  velocity $\bu(\bx) = (-y, x)$, 
pressure $p(\bx) = \lambda x^3 + 0.5(x^2+y^2) -1/3 - \lambda/4$,  force $\bforce(\bx)  = (3 \lambda x^2, 0)$ and $\lambda \geq 0$. We also set the parameter $\lambda = 10$ and the viscosity $\nu = 1$. The velocity gradient is constant throughout the domain.
\item \emph{Polynomial P2P1 test}:  let $(\bu, p) \in [\Pol^2(\Omega)]^2 \times \Pol^{1}(\Omega)$ with velocity $\bu(\bx) = (-x +y^2,y-x^2)$, pressure $p(\bx) = 2x - 2y$ and the corresponding $\bforce(\bx)$ obtained inserting the expression in equation \eqref{eq:Smago:strongform} for $Re = 10000$. 
\end{enumerate}

For clarity in the results analysis, we briefly recall the error estimates from \reviewerB{\cite{BeiraoLovadinaVacca2017}}{\cite{BeiraoLovadinaVacca2018}} for the the norms defined in Section~\ref{sec:NSE}: \reviewerAddB{assuming  $\bu, \bforce \in 
[H^{s+1} (\Omega )]^2$, with $0 < s \leq  k$, then}
\begin{align}
    \left\|\mathbf{u}-\uhv\right\|_{\bbU}  
    &\leq \reviewerB{}{ h^s \mathcal{J}(\bu;\nu) + h^{s+2} \mathcal{H}(\bforce,\nu)} 
    \label{eq:estimator:velocity},\\
    \left\|p-p_h\right\|_{\mathbb{Q}} 
    &\leq \reviewerB{}{C h^s  |p|_{H^{s}(\Omega)} +  C h^{s+2}|\bforce|_{[H^{s+1}(\Omega)]^2}  +  C h^s \mathcal{K}(\bu; \nu),}\label{eq:estimator:pressure}
\end{align}
with \reviewerB{}{the functions $\mathcal{J}$, $\mathcal{H}$, $\mathcal{K}$, and the constant $C$} independent of the meshsize. \reviewerB{}{In particular, $\mathcal{J}$ and $ \mathcal{H}$ depend on $|\bu|_{[H^1(\Omega)]^2}$ and $|\bforce|_{[H^{s+1}(\Omega)]^2}$, respectively.}
\reviewerA{Hence, as long as the force is a polynomial of degree $k$, the discrete velocity is exact.}{In the \emph{Irrotational Force Test}, where both the velocity field and the forcing term are polynomials of degree less than or equal to $k$, we expect the discrete velocity to be reproduced exactly.  
On the other hand, this property is not expected in the \emph{Polynomial P2P1 Test} due to the presence of the Smagorinsky term.} 
\reviewerB{}{Finally, we remark that, in these convergence examples, we employ the estimates for the Navier-Stokes equations proposed in \cite{BeiraoLovadinaVacca2018}, since the Smagorinsky term vanishes for sufficiently small meshsize.} 

We consider uniformly refined mesh sequences. At each recurrence \reviewerA{, a manufactured solution can be obtained from the physical analytical solution and  the Eddy viscosity term}{of the mesh refinement,  the forcing term can be obtained from the analytical solution using the Smagorinsky model}  in its strong form \eqref{eq:Smago:strongform}.
Therefore, for each mesh, we compare the corresponding manufactured analytical solution and the numerical solution. Since we expect the Smagorinsky term to vanish as the mesh size tends to zero due to \eqref{eq:smagoVemterm}, we also expect the error to converge to the physical analytical solution as in \eqref{eq:estimator:velocity}. 
In the irrotational force test, provided that uniform meshes are used\reviewerA{}{, and since the velocity gradient is constant}, $\nuSmagoOper$  remains constant. Since the Laplacian and the gradient of $\nuSmagoOper$ vanish, the force balances only with the pressure gradient and the convective term. This implies that the physical solution and the manufactured problem coincide,  with the latter being independent of $h_T$. 
\begin{figure}[H]
    \centering
    \includegraphics[width=0.5\linewidth]{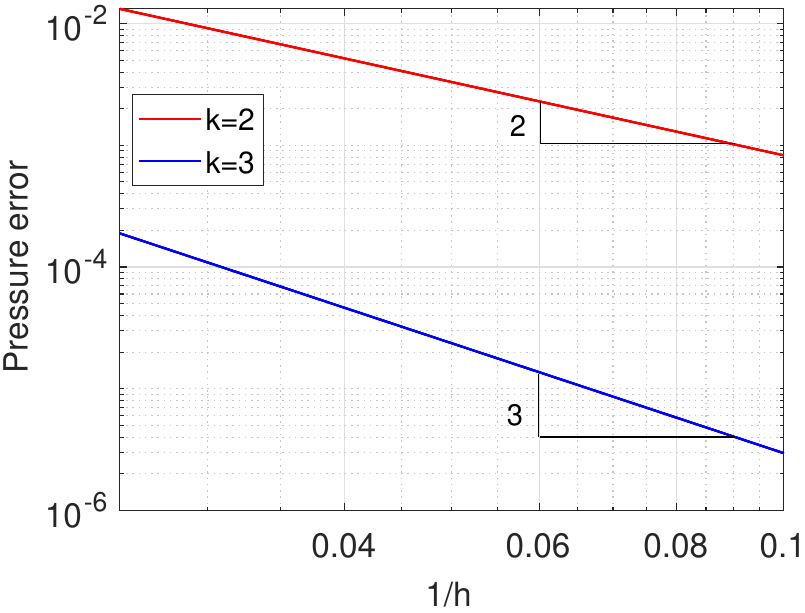}
    \caption{\emph{Irrotational force test}. Pressure errors $\left\| p-p_h\right\|_{\mesh}$ for $Re=1$ and $k \in \{2,3\}$.}
    \label{fig:ExIrrotational:errors}
\end{figure}
In Figure~\ref{fig:ExIrrotational:errors} and 
Table~\ref{tab:ExP2P1:Re10000:errors}, we report the errors and rates of convergence for the irrotational force test and the polynomial P2P1 test cases, respectively. In the irrotational case, we solve using degree $k \in \{2,3,4\}$, but we show pressure errors only for $k=2$ and $k=3$,  as errors for $k=4$ were close to machine precision, as expected. For the same reason, velocity errors for all degrees were omitted. We recall that VEM approximates exactly pressures in $\Pol^{k-1}(\Omega)$,  which for $k=4$ it coincides with the polynomial space of the exact pressure. On the other hand, in the second case, the force $\bforce$ includes non-polynomial components due to the derivative of the Smagorinsky term. This explains the higher velocity errors, although velocity is a polynomial in $[\Pol^2(\Omega)]^2$. This is consistent with \eqref{eq:estimator:velocity} and \eqref{eq:estimator:pressure}.

\begin{table}[h]
    \centering
     \caption{\emph{Polynomial P2P1 test}. Errors and convergence rates using the discrete Smagorinsky model, $Re = 10000$ and $k=2$.}
    \label{tab:ExP2P1:Re10000:errors}
    \resizebox{\textwidth}{!}{%
    \begin{tabular}{ccccccccc}
    \hline
    &    \multicolumn{5}{c}{Velocity}    &&   \multicolumn{2}{c}{Pressure}      \\ \cline{2-6}\cline{8-9}
    $h$     & $\|\nabla \bu - \ProjZeroGrad[k-1]{\uhv}\|_{\mesh}$  &   rate  &&$\| \bu - \ProjZero[k]{\uhv}\|_{\mesh}$ &  rate &&$\| p - p_h\|_{\mesh}$&  rate \\ \cline{2-3} \cline{5-6}\cline{8-9}
    0.100 & 7.3571E-03 & && 4.3530E-04 & && 3.8144E-05 & \\
    0.050 & 1.7860E-03 & 2.043 && 4.1754E-05 & 3.382 && 6.5666E-06 & 2.539 \\
    0.025 & 4.8040E-04 & 1.894 && 5.7269E-06 & 2.866 && 1.7309E-06 & 1.924 \\ \hline
    \end{tabular}
    }
\end{table}

\subsection{Lid-driven cavity flow} For this test case, we focus on a more applied yet established test case in the turbulence model community. The main goal is to provide a first numerical investigation on how VEM approximation performs in a turbulent setting for several Reynolds numbers, highlighting its effectiveness and strengths.

\begin{figure}[h]
    \centering
    \includegraphics[width=0.32\linewidth]{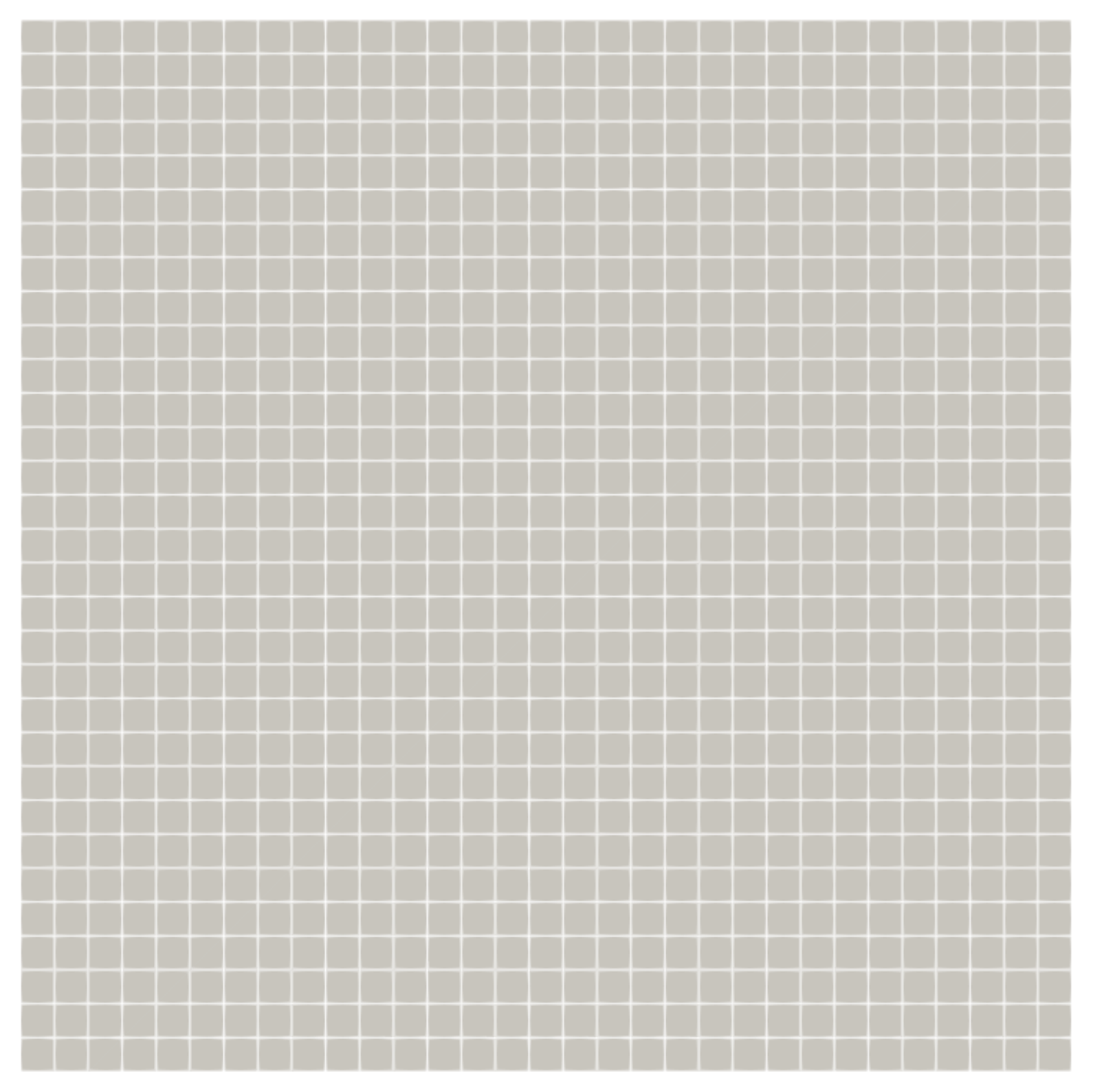}
    \includegraphics[width=0.32\linewidth]{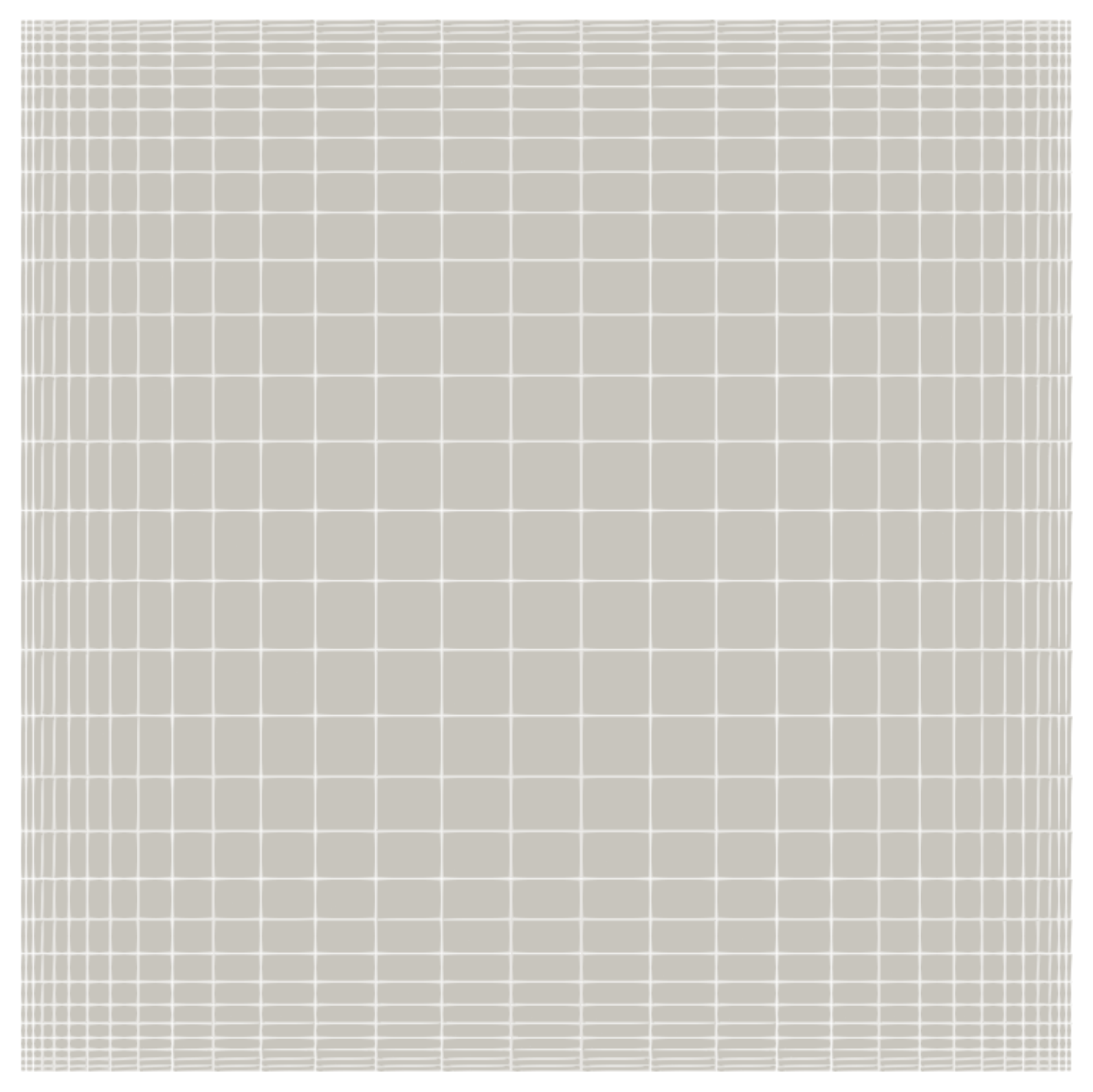}
    \includegraphics[width=0.32\linewidth]{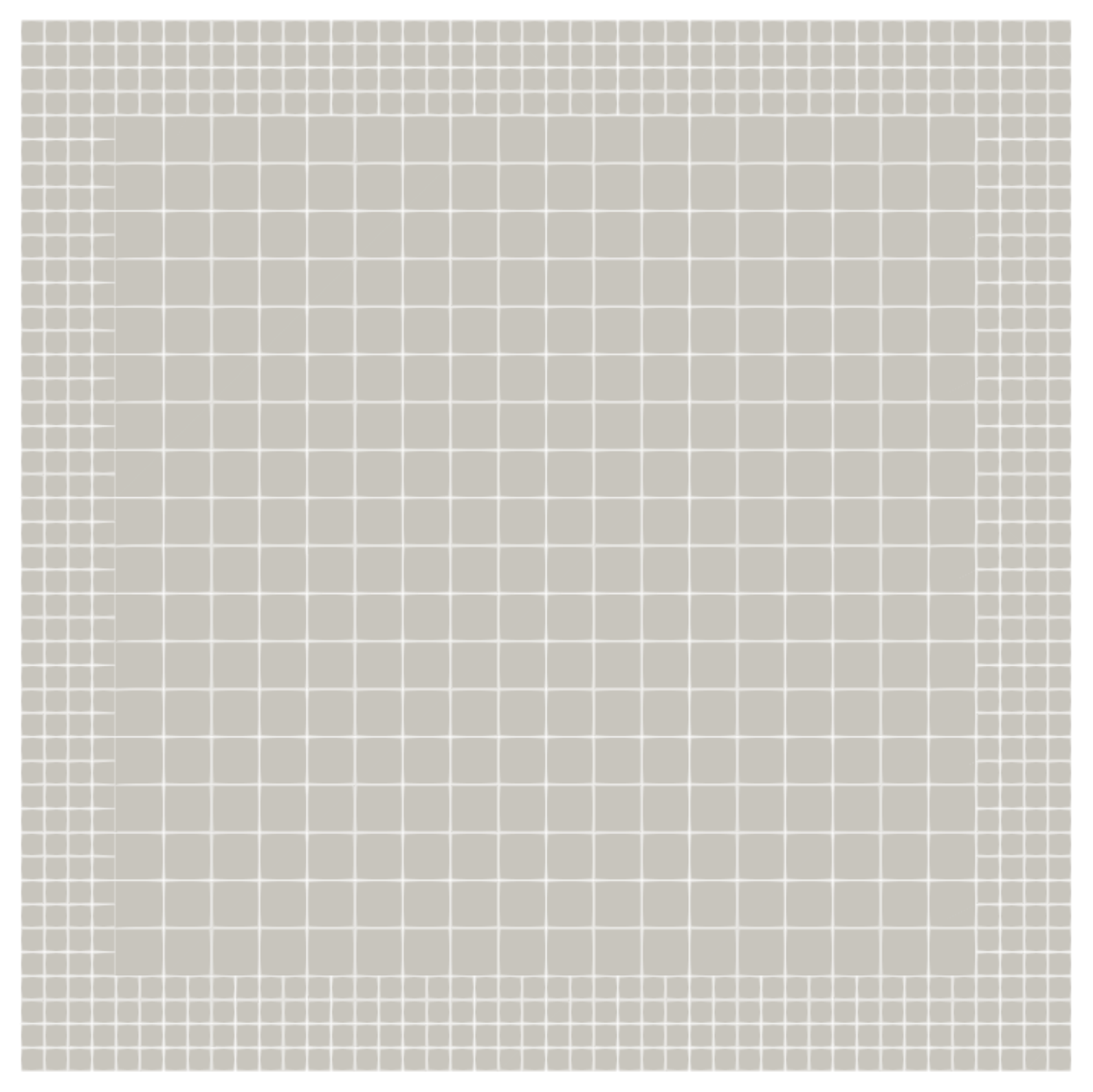}
    \caption{\emph{Lid-driven cavity flow}. Mesh visualization: $31\times31$ USM (left), $31\times31$ ARM  (center), and IMH based on a $22\times22$ USM (right).}
    \label{fig:cavity:meshes}
\end{figure}

In this cavity flow test, no-slip boundary conditions are prescribed in $\Omega$ for the vertical and the bottom walls, while on the top wall a unitary horizontal velocity is imposed. The external body force is set as $\bforce=\bm 0$. The transition regime from laminar to turbulent is ``expected" in the range $6000\leq Re \leq 8000$, see \cite{Erturk2009} and the references therein.
To illustrate the benefits of using VEM, we consider $k=2$ and three different families of meshes: the Uniform Square Meshes (USM), Anisotropic Rectangular Meshes (ARM), and Isotropic Meshes with Hanging nodes (IMH). These families are depicted
in Figure~\ref{fig:cavity:meshes}, while their description at each refinement level is listed in Table $\ref{tab:mesh_setting}$. We remark that the comparison among families is conducted by maintaining a comparable number of DoFs at each refinement level (labeled as Coarse, Medium, and Fine). 
On one hand, the ARMs are obtained by partitioning each side of the domain according to the following function, as in \cite{Haferssas2018}: 
\begin{equation}
    f(\xi)=0.5\left\{1+\frac{\tanh (4\xi-2)}{\tanh (2)}\right\}.
\end{equation}
On the other hand, the IMHs refinement is based on sequences of meshes $(\mathcal{T}_l)_{l>0}$ generated iteratively thanks to the following strategy: let $\mathcal{T}_0$ be a uniform squared mesh generated by a discrete parameter $h_0 \in (0, 1]$.
We generate the coarsest mesh with hanging nodes of the sequence $\mathcal{T}_1$ by splitting by four all the squares whose maximum distance to the boundary in normal direction is $\delta_0 \in [h_0, 0.5 - h_0)$. For our experiments we use $\delta = 0.175$. The consecutive meshes for $\mathcal{T}_l$ are generated using $\delta = \delta_0/2^{l}$. 
The transitions from fine to coarse cells are obtained by using hanging nodes. This approach allows for varying mesh sizes without unnecessarily propagating the refinement throughout the domain, as in the uniform case. We remark that the hanging nodes in the IMHs are easily amenable to VEM, which does not require special treatments. More precisely, the original face is just treated as two collinear faces divided by the hanging node, which is treated as another vertex carrying DoFs.
Here, we mainly use the first level for three different initial meshes, see Table~\ref{tab:mesh_setting} for more details. 

\begin{figure}[h]
    \centering
    \includegraphics[height=6cm, width = 12.5 cm]{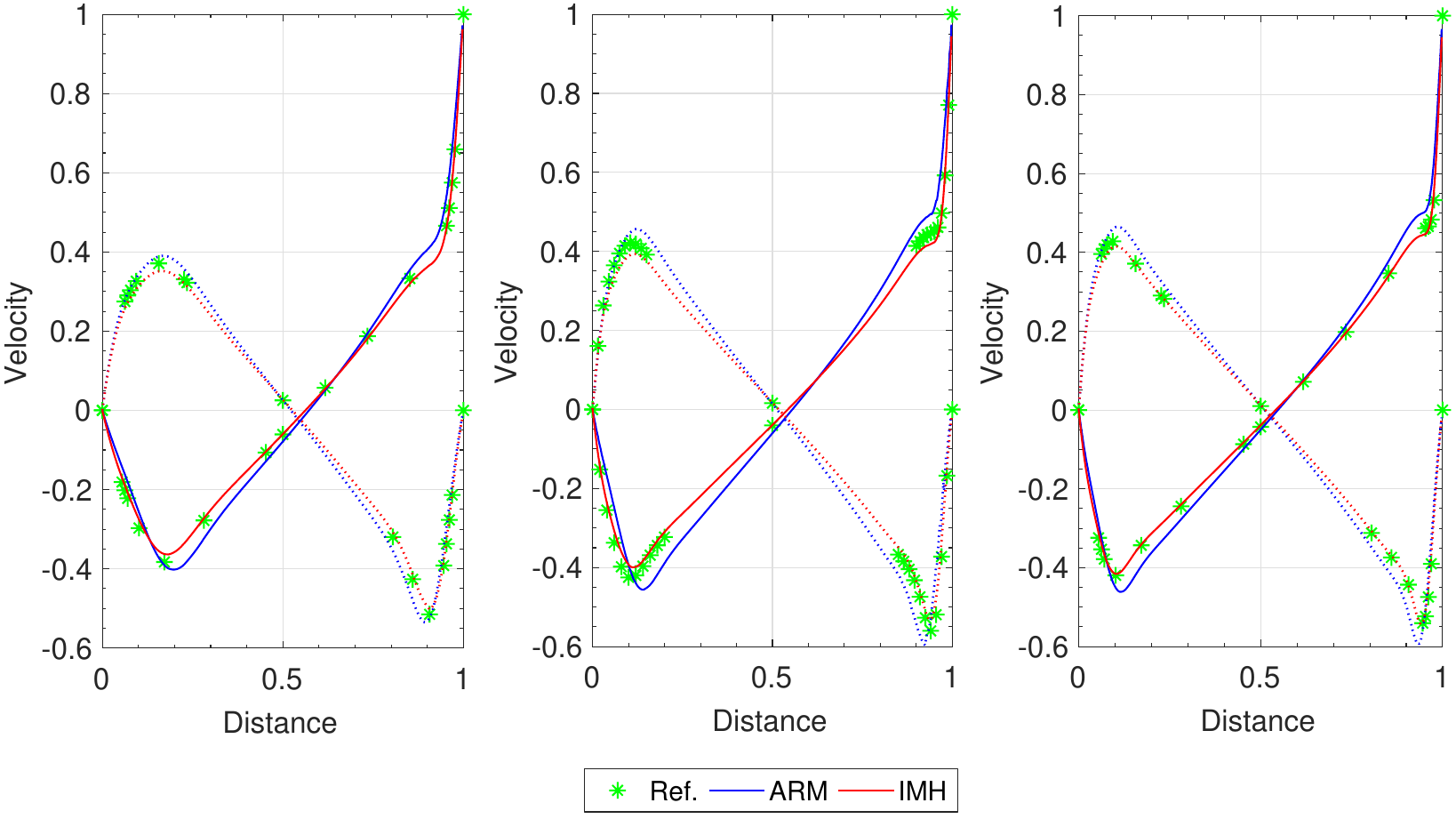}
    \caption{\emph{Lid-driven cavity flow}. Velocity profiles with various Reynolds numbers: $Re = 1000$ using coarse meshes (left), $Re=2500$ using medium meshes (center), and $Re=3200$ using fine meshes (right). The solid and dotted lines represent the $x$- and $ y$-components of the velocity, respectively.}
    \label{fig:cavity:profiles:anisotropy}
\end{figure}

In the literature, the ARMs are commonly used to achieve higher accuracy near the domain walls, where high gradients and small scales are expected as the Reynolds number increases. However, as one can observe from the middle plot of Figure~\ref{fig:cavity:meshes}, some elements are much stretched, and this anisotropy of the elements near the boundary of the domain might foster the overestimation of the Smagorinsky viscosity near the walls. Indeed, the Smagorinsky model can be defined as \emph{isotropic}: the model's scaling $\hCell$ does not take into account the element anisotropy, treating every direction of the cell similarly.
In this specific boundary layer, the cell diameter $\hCell$ of each element of the layer is, roughly speaking, comparable to the measure of its largest face. 
This phenomenon is visible in Figure~\ref{fig:cavity:profiles:anisotropy}, where we plot the $x$- and $ y$-~components of the velocity. The ARM solution (blue solid and blue dotted lines) does not match the Erturk et al. reference values (green stars) \cite{Erturk2009} for $Re=1000$ (coarse mesh), $Re=2500$ (medium mesh), and $Re=3200$ (fine mesh). We show the highest Reynolds number that each mesh can resolve, and for which reference data is available for comparison.
We conclude that ARMs are not practical in this setting. Alternatively, uniform meshes, as USM, would require a significantly high number of elements to adequately resolve localized features such as sharp gradients and boundary layers, leading to unnecessary computational cost. Instead, a non-uniform mesh allows for targeted refinement in critical areas, improving accuracy while maintaining overall efficiency.
Thus, to assess the turbulent model without the influence of the mesh anisotropy and focusing on the boundary layer refinement, we use the IMH family of meshes.

\begin{table}[h]
    \centering 
    \caption{\emph{Lid-driven cavity flow}. Mesh setting description, where $N$ denotes the number of segments into which the domain is partitioned in each direction.}
    \begin{tabular}
    {lcccc}
    \hline
    & \multicolumn{2}{c}{IMH}& \multicolumn{2}{c}{ USM/ARM}\\
    & $N\times N$ & $\#$ DoFs & $N\times N$ &$\#$ DoFs\\ \hline
     Coarse & 22 $\times$ 22 & 10399 & 31 $\times$ 31 & 10526 \\
     Medium & 28 $\times$ 28 & 15235 & 38 $\times$ 38 & 15583 \\
    Fine & 40 $\times$ 40 & 31887 & 54 $\times$ 54 & 31647 \\ \hline
    \end{tabular}
    \label{tab:mesh_setting}
\end{table}

\begin{figure}[htp]
    \centering
    \begin{subfigure}[b]{0.475\textwidth}
    \includegraphics[width=.9\linewidth]{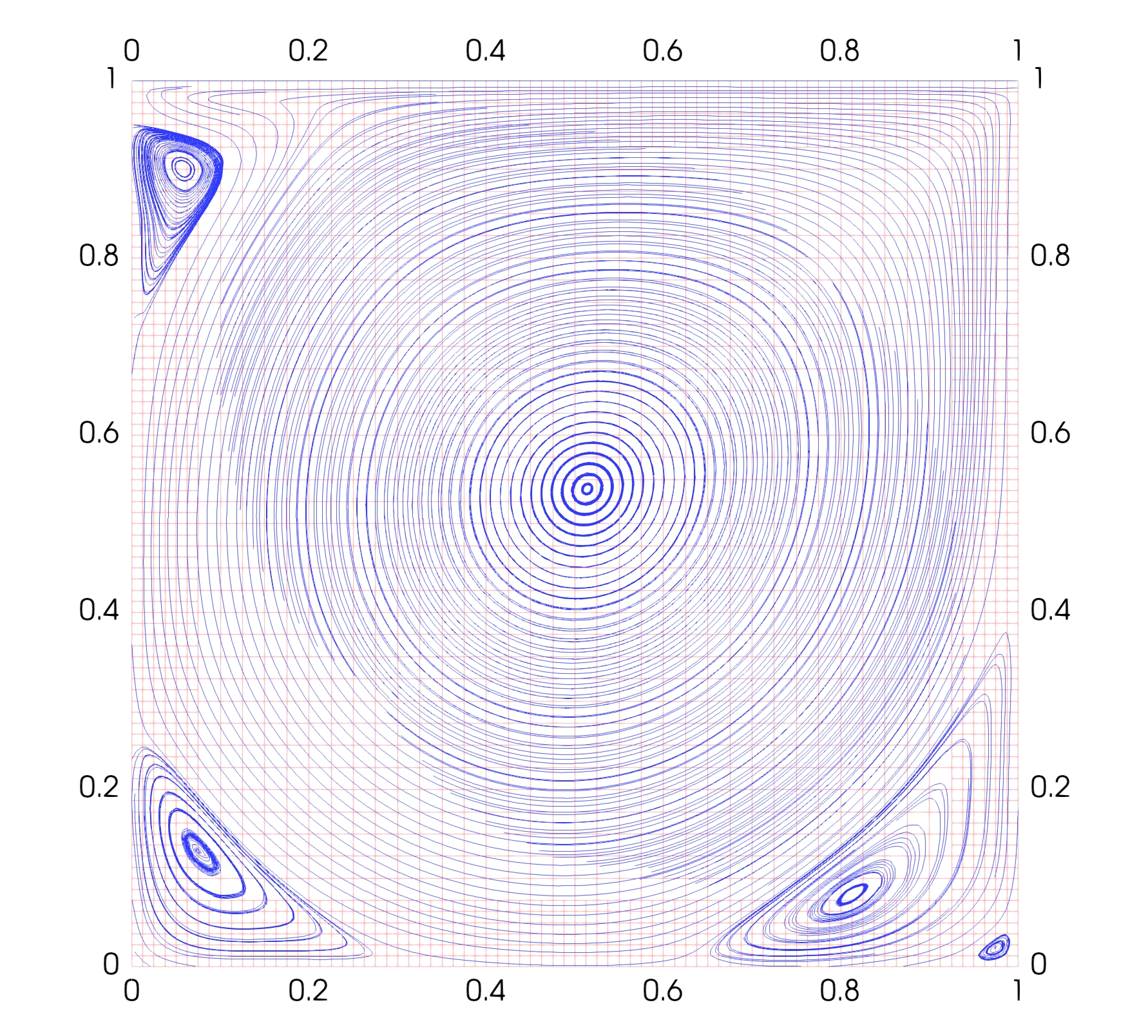}
       \caption{$Re = 5000$ and fine IMH \\ \reviewerAddA{(Smagorinsky)}}
    \end{subfigure}
   \begin{subfigure}[b]{0.475\textwidth}
   \includegraphics[width=.9\linewidth]{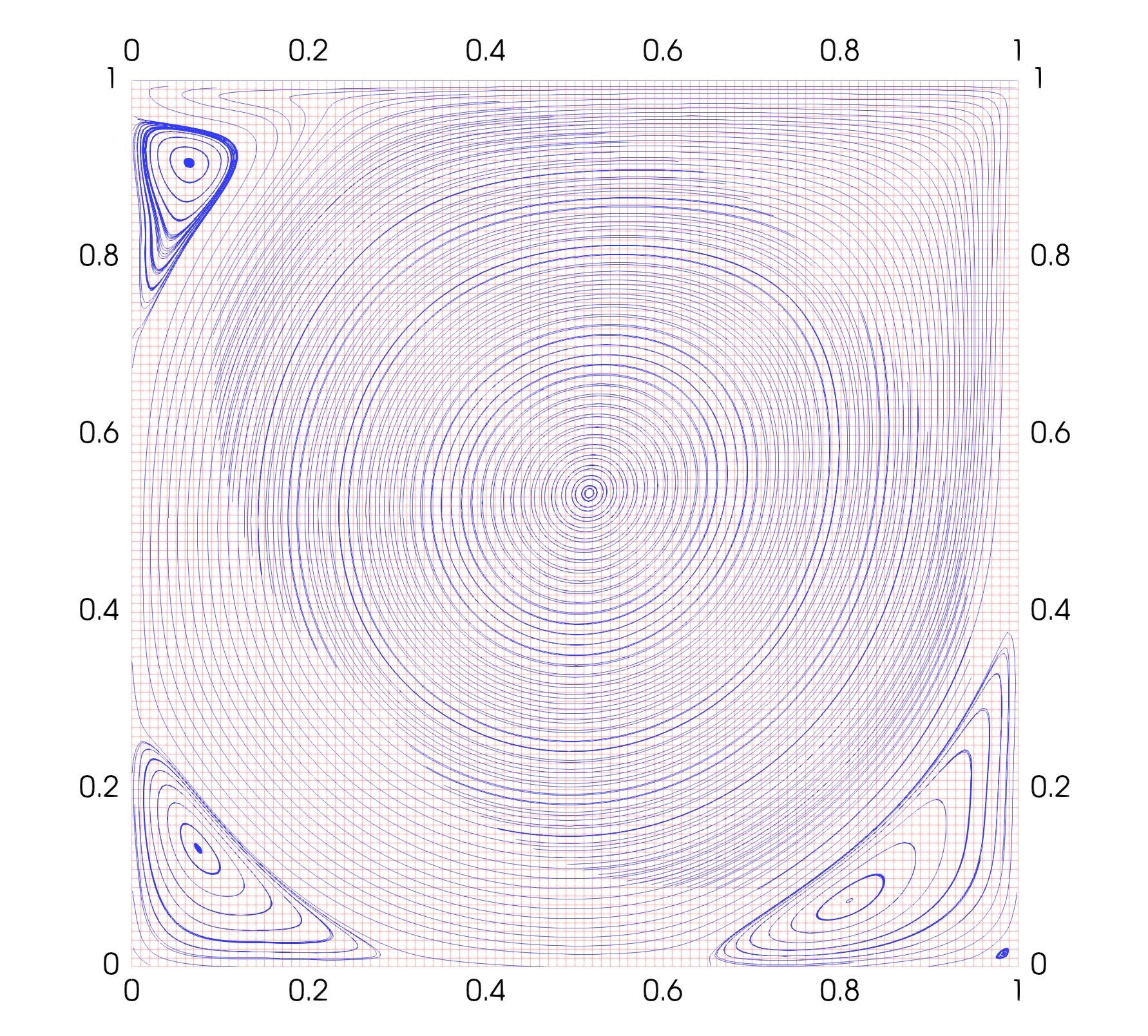}  
   \caption{$Re = 5000$ and USM 100x100 \\\reviewerAddA{(NSE)}}
    \end{subfigure}

    \begin{subfigure}[b]{0.475\textwidth}
    \includegraphics[width=.9\linewidth]{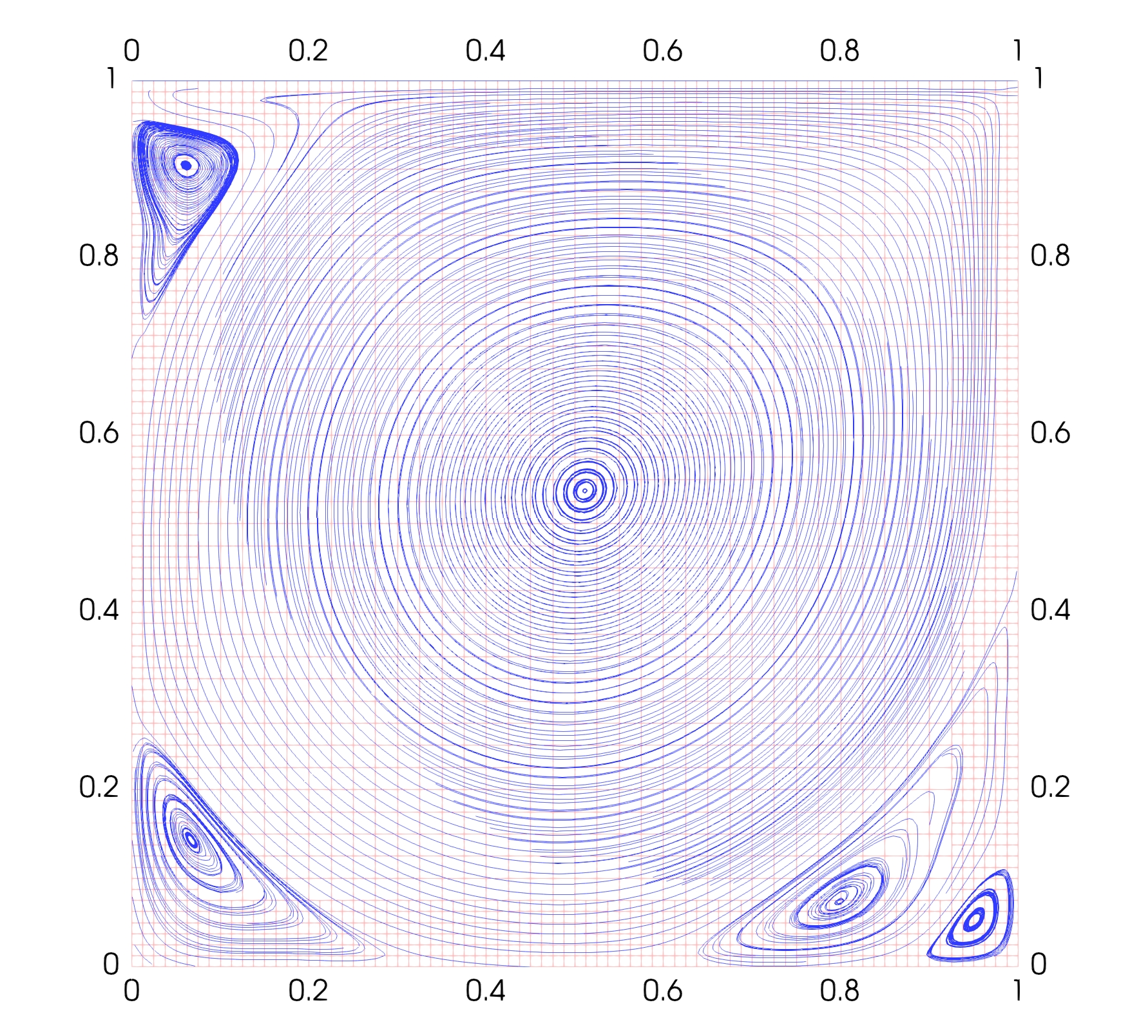}
       \caption{$Re = 7500$ and fine IMH \\\reviewerAddA{(Smagorinsky)}}
    \end{subfigure}
   \begin{subfigure}[b]{0.475\textwidth}
   \includegraphics[width=.9\linewidth]{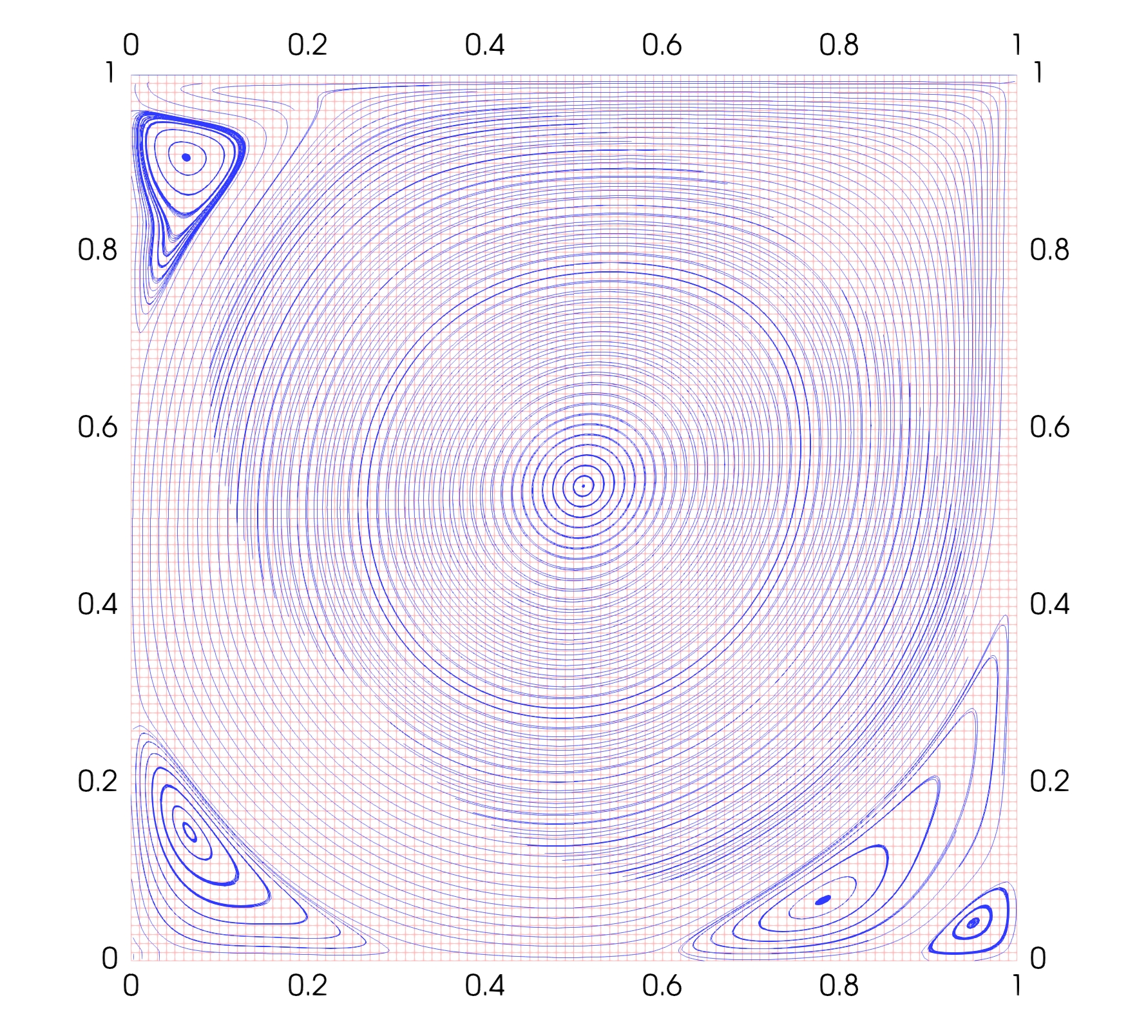}  
   \caption{$Re = 7500$ and USM 100x100 \\ \reviewerAddA{(NSE)}}
    \end{subfigure}
    \begin{subfigure}[b]{0.475\textwidth}
    \includegraphics[width=.9\linewidth]{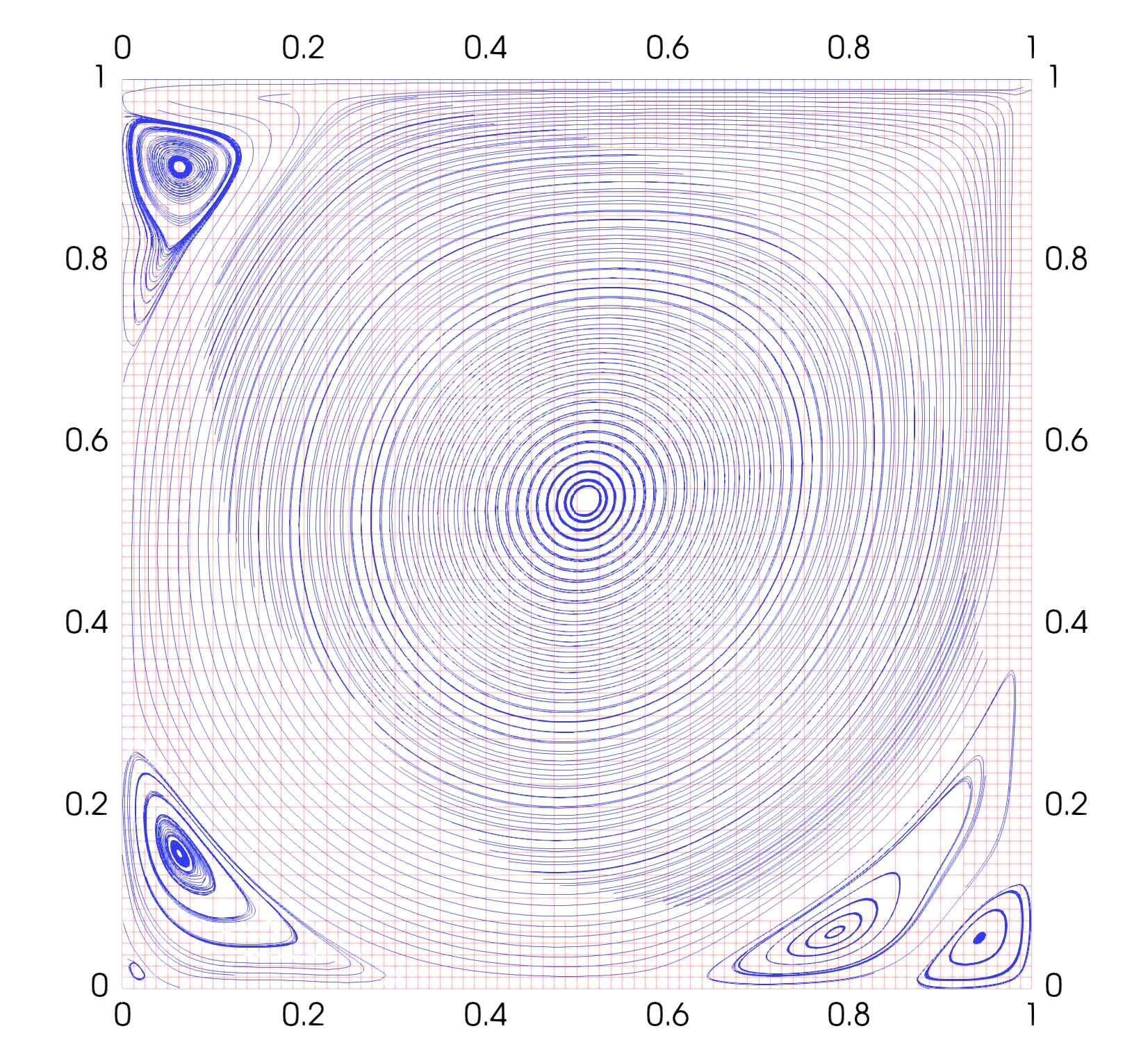}
    \caption{$Re = 10000$ and fine IMH \\\reviewerAddA{(Smagorinsky)}}
     \end{subfigure}
     \begin{subfigure}[b]{0.475\textwidth}\includegraphics[width=.9\linewidth]{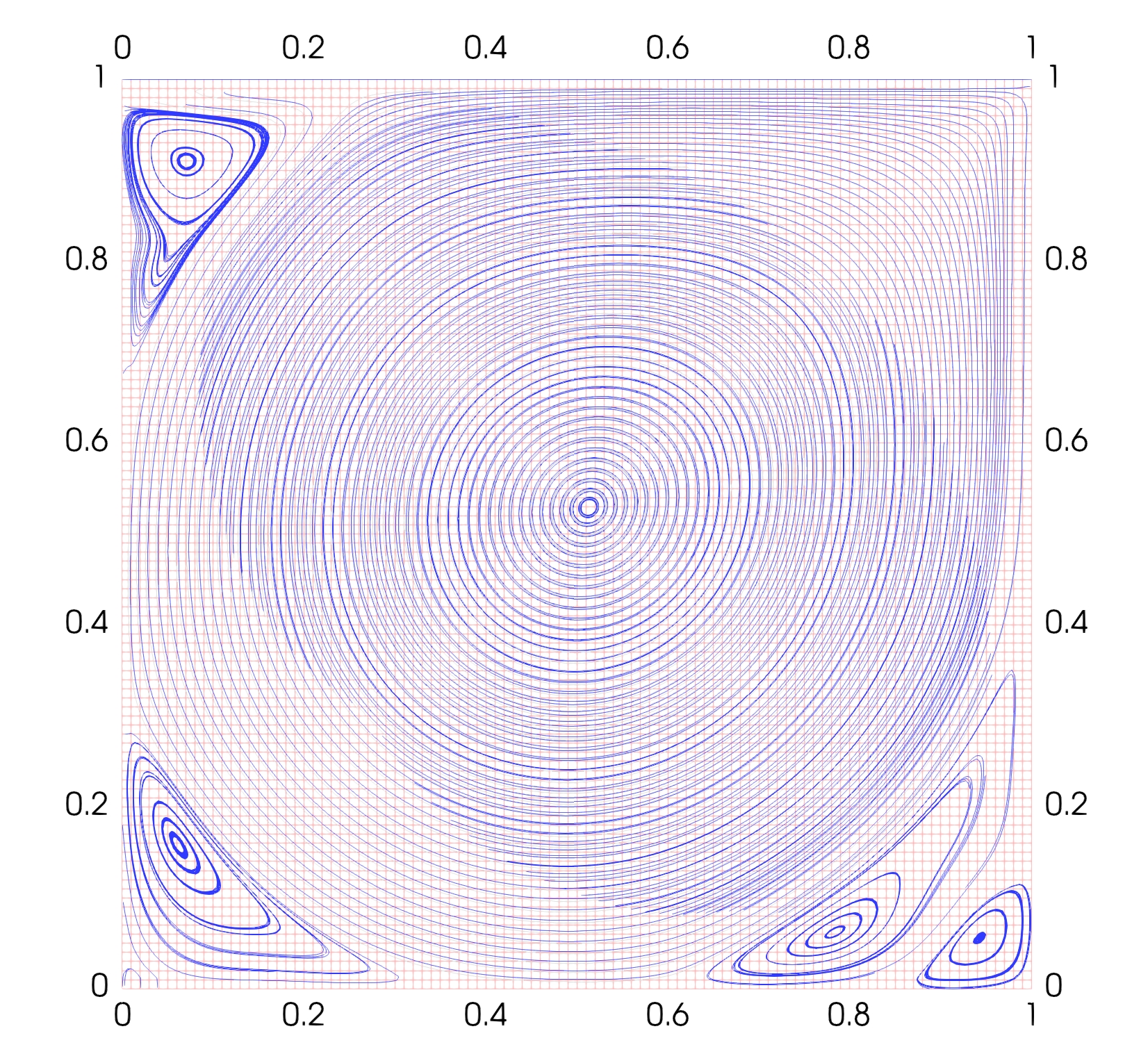}
     \caption{$Re = 10000$ and USM 100x100 \\\reviewerAddA{(NSE)}}
     \end{subfigure}
     \caption{\emph{Lid-driven cavity flow}. Streamlines for several Reynolds numbers, using the fine IMH mesh \reviewerAddA{ with Smagorinky model} compared to the USM 100x100 \reviewerAddA{with NSE}. }
    \label{fig:cavity:streamlines}
\end{figure}

From Figure~\ref{fig:cavity:profiles:anisotropy}, we observe that IMH results (red solid and red dotted lines) better match the reference values with respect to ARM results. We do not report the results for USM since they were comparable to IMH in terms of accuracy. This is expected since the uniform mesh has a good description of the boundary layer, but increases the number of DoFs, more than twice with respect to the finer mesh, i.e., IMH presents 31887 DoFs, while USM presents 69763 DoFs. This is caused by a non-necessary refinement on the central region of the domain.
\begin{figure}[h]
    \centering
    \begin{subfigure}[b]{0.475\textwidth}
    \includegraphics[width=.9\linewidth]{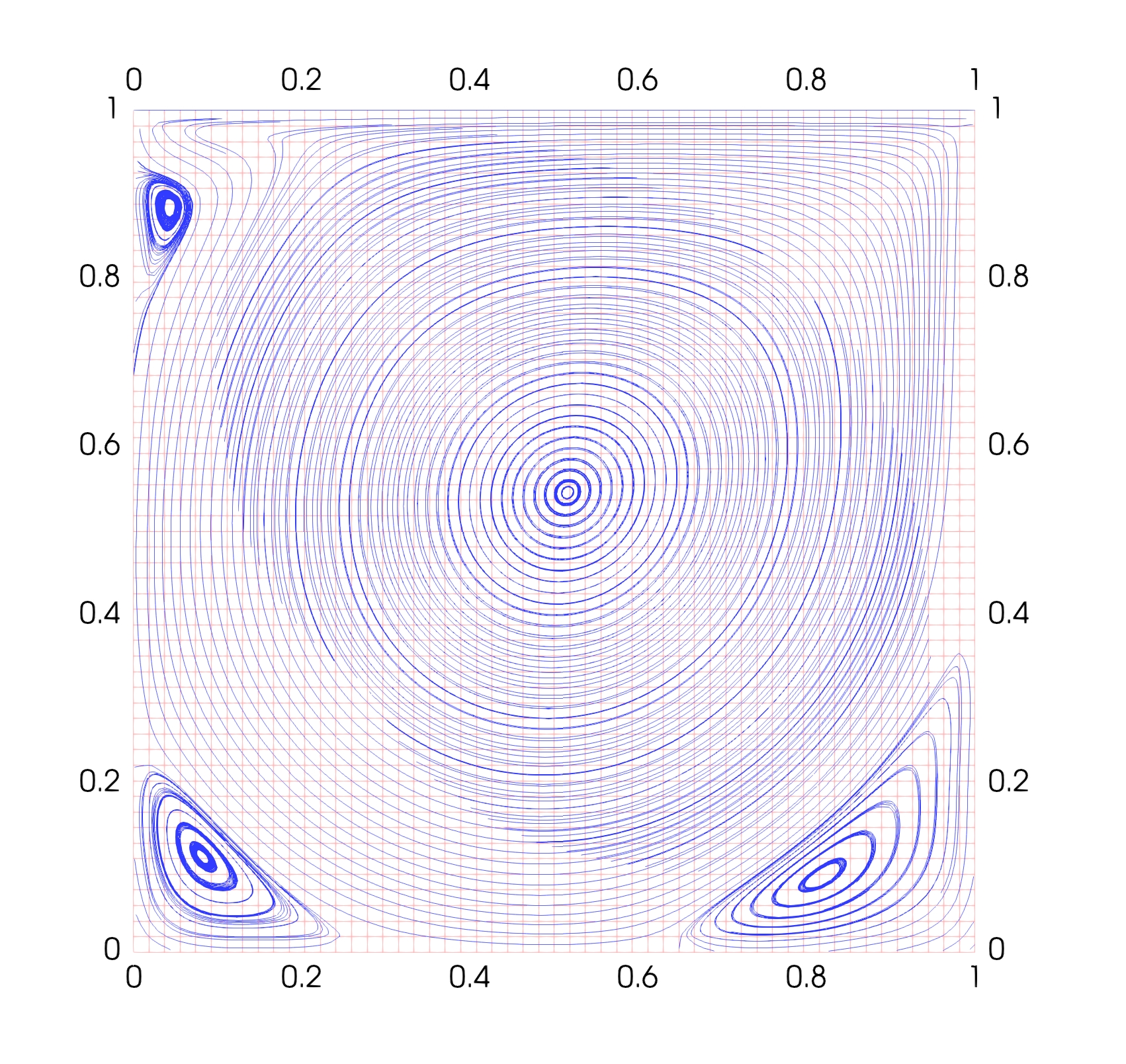
    }
       \caption{Fine USM \reviewerA{}{(NSE)}}
    \end{subfigure}
   \begin{subfigure}[b]{0.475\textwidth}
   \includegraphics[width=.9\linewidth]{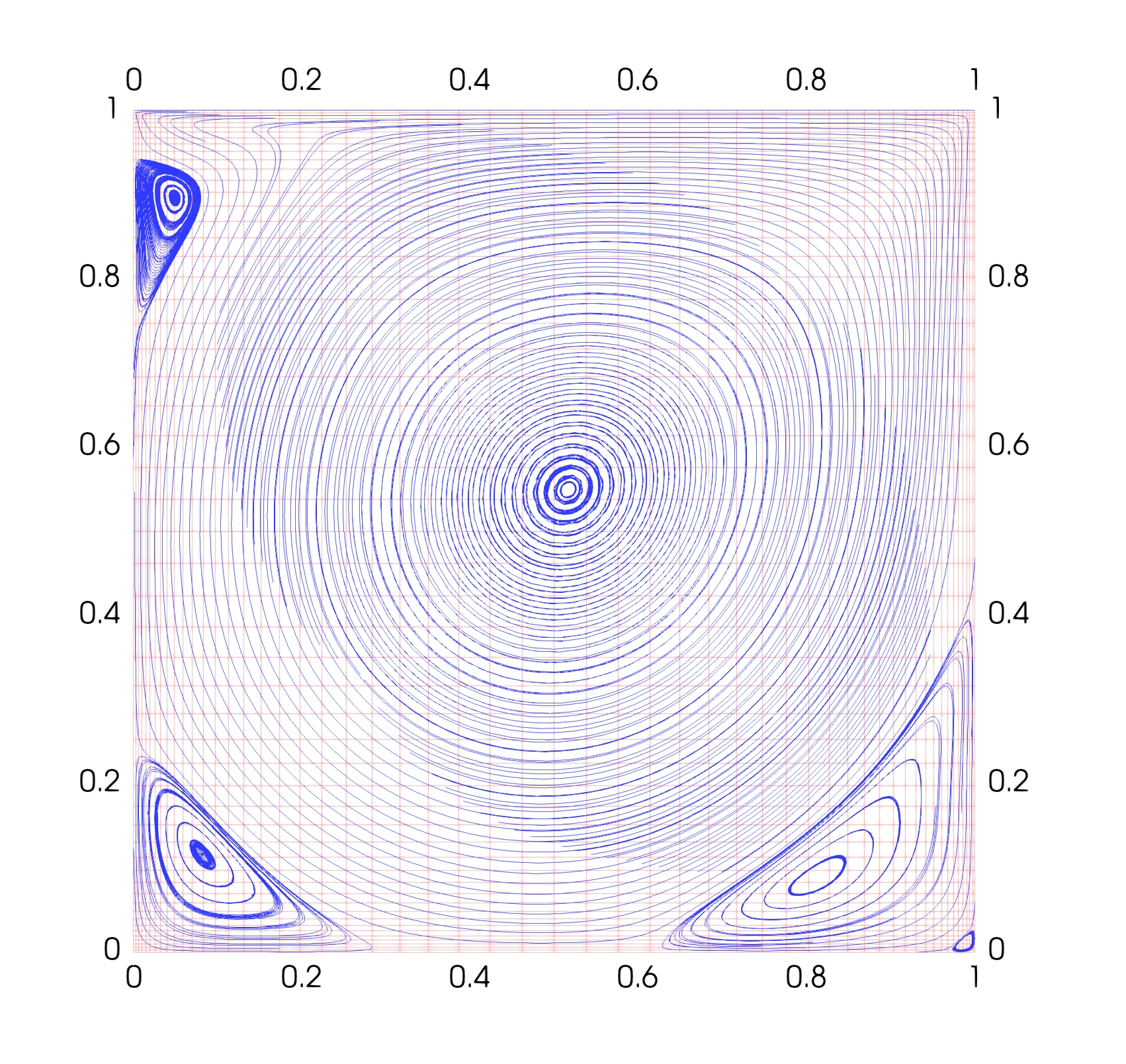}  
   \caption{Fine ARM \reviewerA{}{(Smagorinsky)}}
    \end{subfigure}
    \begin{subfigure}[b]{0.475\textwidth}
    \includegraphics[width=.9\linewidth]{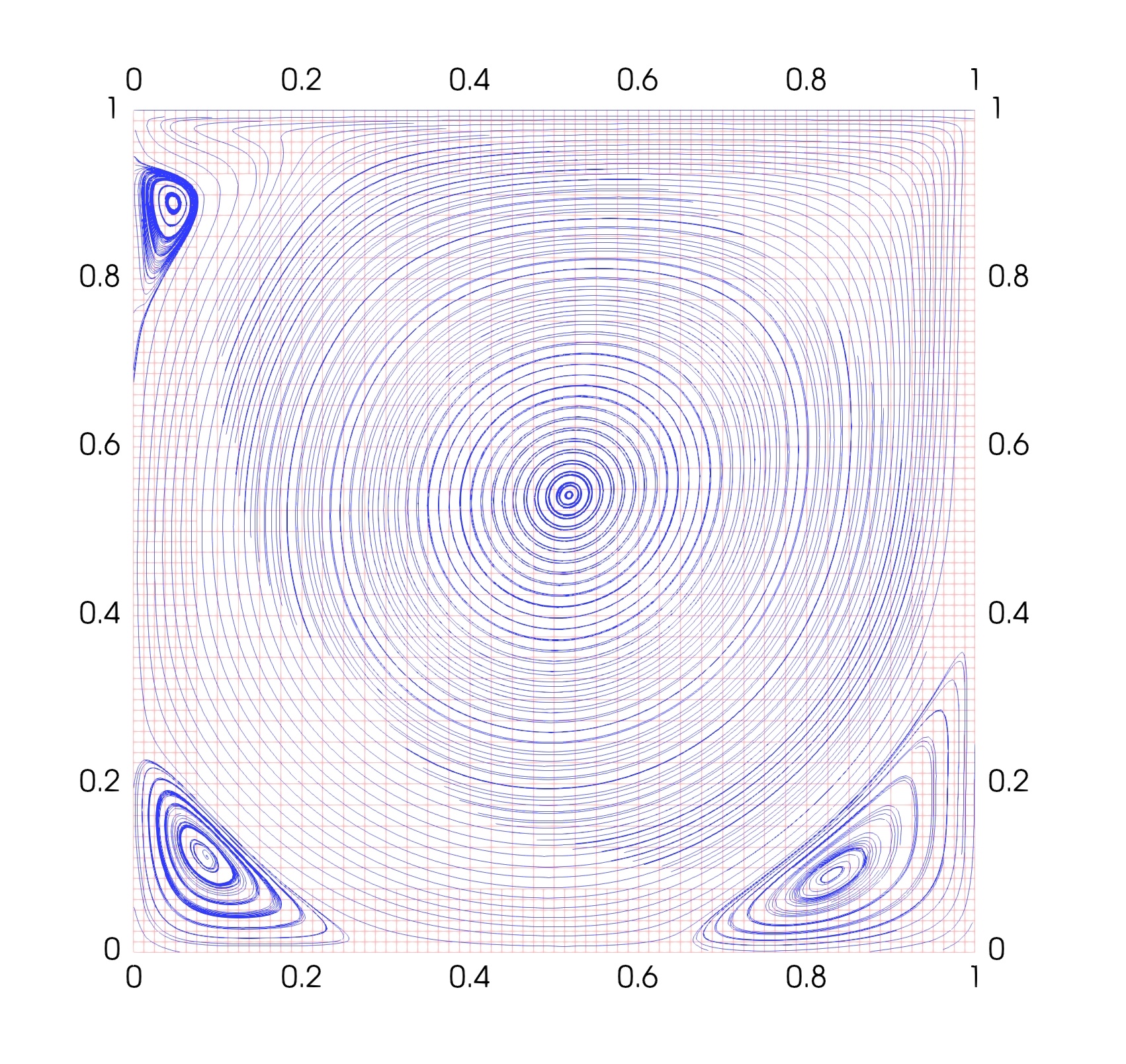}
    \caption{Fine IMH \reviewerA{}{(Smagorinsky)}}
     \end{subfigure}
     \begin{subfigure}[b]{0.475\textwidth}\includegraphics[width=.9\linewidth]{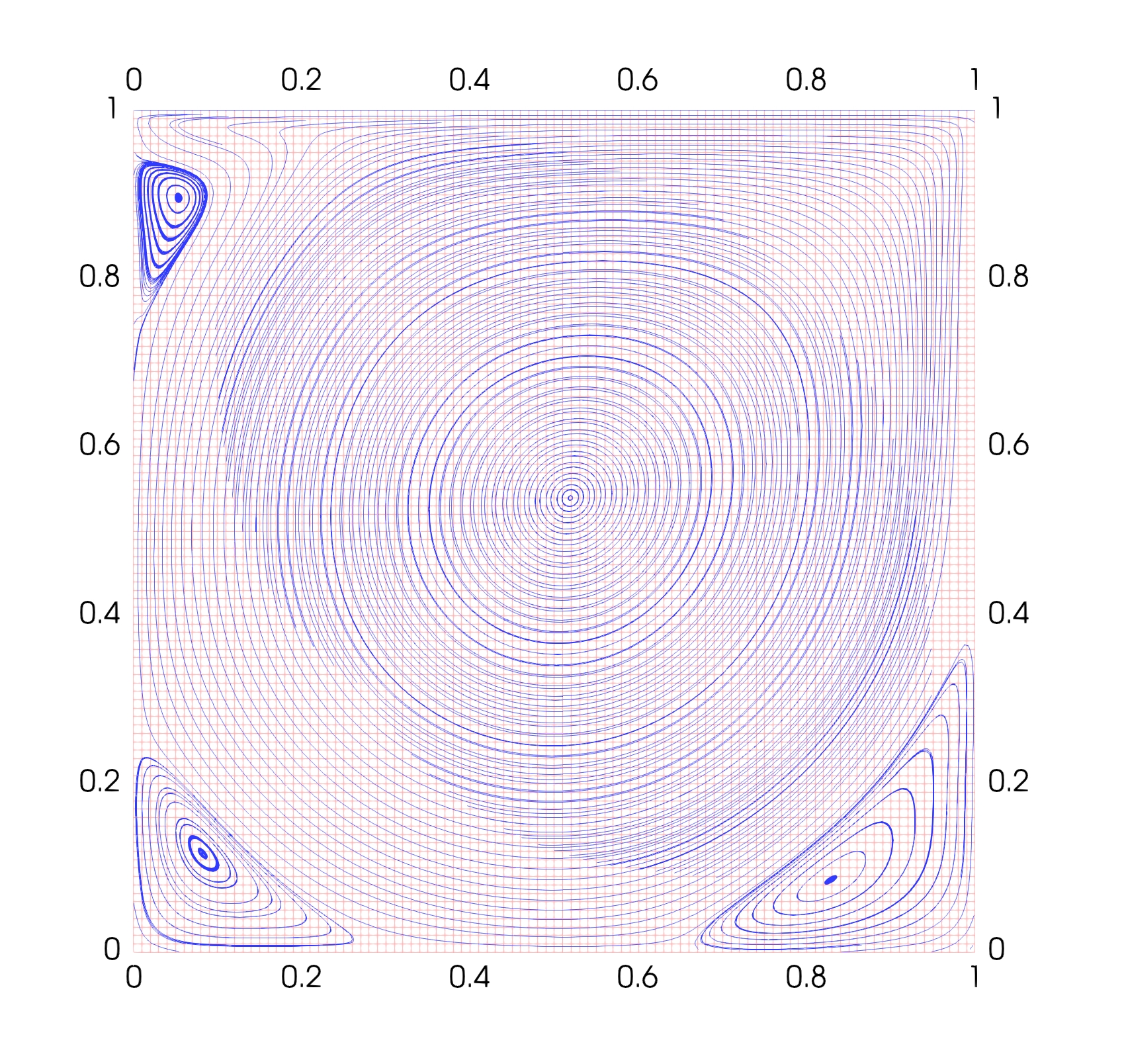}
     \caption{USM 100x100 \reviewerA{}{(Smagorinsky)}}
     \end{subfigure}
     \caption{\emph{Lid-driven cavity flow}. Streamlines for  $Re=3200$ \reviewerAddA{with the Smagorinsky model} using the three fine meshes, and the USM 100x100 \reviewerAddA{with NSE} for comparison \reviewerAddA{(top left)}.}
    \label{fig:cavity:streamlines2}
\end{figure}
{In Figure~\ref{fig:cavity:streamlines2}, we show the streamlines for a laminar case with $Re=3200$ and the finest mesh for each mesh family. We qualitatively compare the results and observe that in the case of ARM, it shows the biggest displacements of the vortices with respect to the reference, which are twice and even three times the displacements shown by the IMH or USM. In addition, there is a spurious vortex in the bottom-right corner. In Figure~\ref{fig:cavity:streamlines}, we show the streamlines plot for $Re \in \{5000, 7500, 10000\}$, using the fine IMH \reviewerA{}{for the Smagorinsky model} in comparison with reference solutions based on the standard Navier-Stokes problem solved on a USM $100 \times 100$ with 109203 DoFs. We observe that the Smagorinsky model solution accurately represents the feature of the turbulent flow, and we only see a slight displacement of the vortices of the order of the cell diameter of the reference mesh.}
{To quantitatively assess the convenience of exploiting IMH in this setting, we compute the $L^2$-errors of the velocity with respect to the reference configuration for $Re \in \{5000, 7500, 10000\}$. In Tables \ref{tab:cavity:L2_errors:USM} and \ref{tab:cavity:L2_errors:IMH}, we report the relative errors for the Smagorinsky model applied to the USM and IMH  families, respectively. Henceforth, the hyphen symbol ``-" in the tables indicates a failure of the Newton solver to converge. The IMH family presents smaller errors compared to the USM family and is more robust in terms of the Newton solver convergence: we reach convergence for all the Reynolds numbers and meshes considered, except for $Re=10000$ with the coarse mesh. In contrast, the USM family encounters more convergence issues ($Re \in \{7500, 10000\}$ with the coarse mesh and $Re=10000$ with the medium mesh). We will further comment on the convergence of the Newton solver later in this section.}

\begin{figure}[h]
    \centering
    \includegraphics[height=6cm, width = 12 cm]{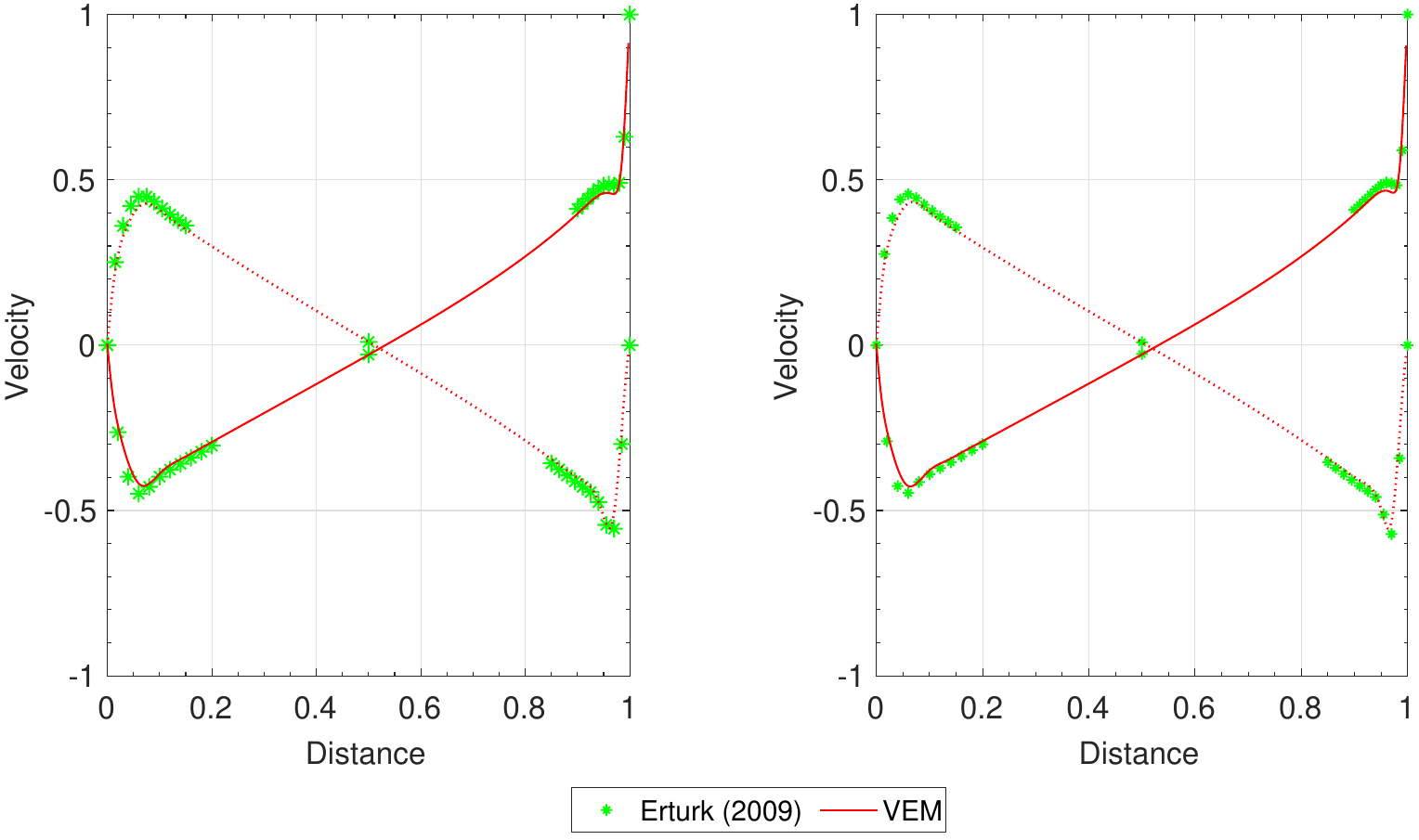}
    \caption{\emph{Lid-driven cavity flow}. Velocity profiles using the fine IMH for $Re=7500$ (left) and $Re=10000$ (right). The solid and dotted lines represent the $x$- and $ y$-components of the velocity, respectively.}
    \label{fig:cavity:velocity_profiles:turb}
\end{figure}

In addition,
in Figure~\ref{fig:cavity:velocity_profiles:turb}, we show the velocity profiles along the vertical and horizontal middle axes for laminar and turbulent (including transitional) flow, respectively. The numerical results are compared with those reported by Erturk et al. \cite{Erturk2009}. We prefer it compared to Ghia et al. \cite{GhiaGhiaShin1982} since it privileges data points close to walls, where sharp gradients are especially steep, and gives more relevant information about how well the refinement captures them.

\begin{table}[]
\centering 
    
    \caption{\emph{Lid-driven cavity flow}. Relative $L^2$-errors \reviewerAddA{for the velocity variable} against the numerical reference solution varying the $Re$ and considering USM. }
    \begin{tabular}{crrr}
    \hline
    &\multicolumn{3}{c}{$Re$}\\ \cline{2-4}
    Mesh& \multicolumn{1}{c}{5000} & \multicolumn{1}{c}{7500}  & \multicolumn{1}{c}{10000}\\ \hline 
    Coarse &  3.071E-01 &           - &           - \\
    Medium &  2.490E-01 &  3.264E-01 &           - \\
    Fine   &  1.809E-01 &  2.433E-01 &  2.941E-01 \\
   \hline
    \end{tabular}
    \label{tab:cavity:L2_errors:USM}
\end{table}
    
\begin{table}[]
    \caption{\emph{Lid-driven cavity flow}. Relative $L^2$-errors \reviewerAddA{for the velocity variable} against the numerical reference solution varying the $Re$ and considering IMH. }
\centering 
    
    \begin{tabular}{crrr}
    \hline
    &\multicolumn{3}{c}{$Re$}\\ \cline{2-4}
    Mesh& \multicolumn{1}{c}{5000} & \multicolumn{1}{c}{7500}  & \multicolumn{1}{c}{10000}\\ \hline 
    Coarse &  1.853E-01 &  2.399E-01 &          - \\
    Medium &  1.449E-01 &  1.942E-01 &  2.351E-01 \\
    Fine   &  1.048E-01 &  1.495E-01 &  1.868E-01 \\
   \hline
    \end{tabular}
    \label{tab:cavity:L2_errors:IMH}
\end{table}

We now propose a quantitative analysis based on the relative $\ell^2$-error against the Erturk et al. \cite{Erturk2009} values for USM meshes, IMH meshes, and ARMs with  $Re \in \{1000, 2500, 5000, 7500, 10000\}$. The results are collected in Tables \ref{tab:cavity:l2_errors:USM}, \ref{tab:cavity:l2_errors:IMH}, and \ref{tab:cavity:l2_errors:newARM}, respectively. Notice that the precision of the error is limited by the precision of the reported reference data in \cite{Erturk2009,GhiaGhiaShin1982}, i.e., five significant digits. 
It is clear that the IMH family offer the best performance in terms of accuracy, providing smaller or, in the worst case, comparable errors with respect to USM or ARM, when comparison is possible, i.e., only when the Newton solver converges for the two considered different meshes.
\begin{table}[h]
    \centering
    \caption{\emph{Lid-driven cavity flow}. Relative $\ell^2$-errors \reviewerAddA{for the velocity variable} against Erturk et al.\ \cite{Erturk2009} varying the $Re$ and considering USM. }
    \begin{tabular}{crrrrr}
    \hline
    &\multicolumn{5}{c}{$Re$}\\ \cline{2-6}
    Mesh&  \multicolumn{1}{c}{1000} & \multicolumn{1}{c}{2500}&  \multicolumn{1}{c}{5000} & \multicolumn{1}{c}{7500}  & \multicolumn{1}{c}{10000}\\ \hline 
    Coarse&  9.365E-02 & 1.538E-01  & 2.399E-01 & -          & -\\
    Medium&	 7.502E-02 & 1.195E-01  & 1.859E-01 & 2.444E-01  & -\\
    Fine  &  5.303E-02 & 7.594E-02  & 1.229E-01 & 1.638E-01  &  1.993E-01 \\
   \hline
    \end{tabular}
    \label{tab:cavity:l2_errors:USM}
\end{table}

\begin{table}[h]
    \centering
    \caption{\emph{Lid-driven cavity flow}. Relative $\ell^2$-errors \reviewerAddA{for the velocity variable} against Erturk et al.\ \cite{Erturk2009} varying the $Re$ and considering IMH. }
    \begin{tabular}{crrrrr}
    \hline
    &\multicolumn{5}{c}{$Re$}\\ \cline{2-6}
    Mesh&  \multicolumn{1}{c}{1000} & \multicolumn{1}{c}{2500}&  \multicolumn{1}{c}{5000} & \multicolumn{1}{c}{7500}  & \multicolumn{1}{c}{10000}\\ \hline 
    Coarse & 7.307E-02 & 1.018E-01  & 1.445E-01 &  1.782E-01 & 1.104E-01 \\
    Medium & 5.913E-02 & 7.813E-02  & 1.104E-01 &  1.385E-01 & 1.626E-01\\
    Fine   & 4.173E-02 &  5.194E-02 & 7.070E-02 &  9.266E-02 & 1.133E-01\\
   \hline
    \end{tabular}
    \label{tab:cavity:l2_errors:IMH}
\end{table}

\begin{table}[h]
    \centering
    \caption{\emph{Lid-driven cavity flow}. Relative $\ell^2$-errors \reviewerAddA{for the velocity variable} against Erturk et al.\ \cite{Erturk2009} varying the $Re$ and considering ARM. }
    \begin{tabular}{crrrrr}
    \hline
    &\multicolumn{5}{c}{$Re$}\\ \cline{2-6}
    Mesh&  \multicolumn{1}{c}{1000} & \multicolumn{1}{c}{2500}&  \multicolumn{1}{c}{5000} & \multicolumn{1}{c}{7500}  & \multicolumn{1}{c}{10000}\\ \hline 
    Coarse & 9.429E-02 & -  & - &  - & - \\
    Medium & 6.972E-02 & 1.602E-02  & - &  - & - \\
    Fine   & 1.054E-02 & 1.054E-02 & - &  - & - \\
   \hline
    \end{tabular}
    \label{tab:cavity:l2_errors:newARM}
\end{table}

\begin{figure}[H]
     \centering
     \includegraphics[width=0.7\linewidth]{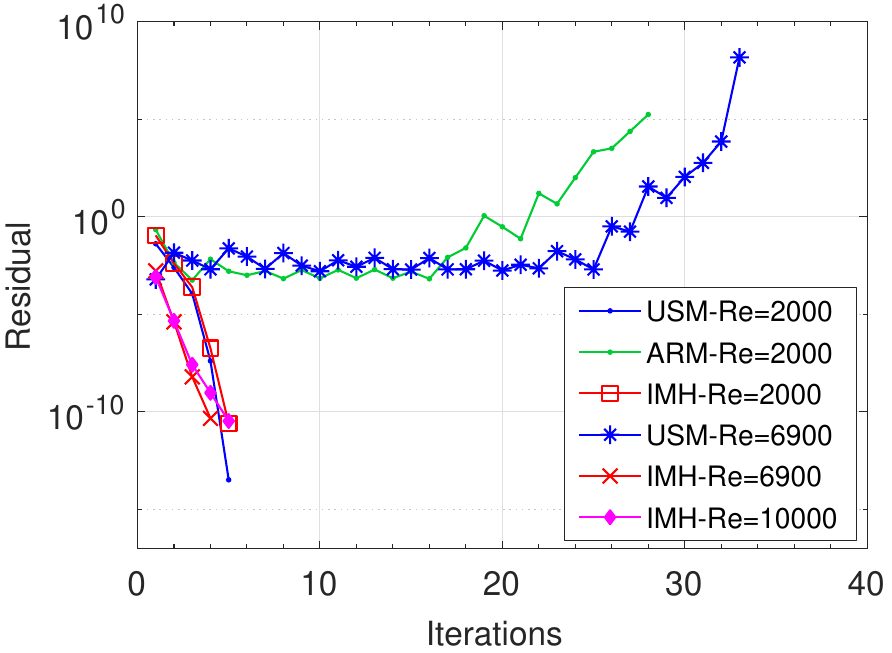} 
     \caption{\emph{Lid-driven cavity flow}. The Newton solver convergence for various $Re$ and the three mesh families.}
     \label{fig:cavity:convergence:coarse}
 \end{figure}

Specifically, we solve the nonlinear system using the Newton solver without preconditioning or relaxation techniques. The Newton solver is well known for exhibiting quadratic convergence, provided that a sufficiently good initialization is supplied. We noticed that, under the same initial guess condition for the solver, the IMH family are the most robust choice and allow for simulations for large Reynolds, up to $Re=10000$. On the contrary, ARMs family yield convergence issues. 
In Figure~\ref{fig:cavity:convergence:coarse}, we show the convergence process for the three coarse meshes for different representative Reynolds. From the plot, we observe that for $Re=2000$, IMH and USM converge in the same number of iterations, while the ARM already diverges. When we increase the Reynolds number to $Re=6900$, the USM diverges. The IMH converges in few iterations not only for $Re=6900$ but also for large Reynolds, up to $Re=10000$. 

\begin{remark}[Anisotropic scaling]
Let $T \in \mesh$ be an element of the mesh. We recall that $h^\star_{\face}\coloneqq \displaystyle \min_{\face \subset \dCell}h_F$ denotes the smallest face diameter of $T$. In the numerical investigation we carried out, choosing $h_T$ in the Smagorinsky model for ARMs led to overdiffusive results due to the anisotropic nature of the mesh. For this reason, we tried $h^\star_{\face}$ as an alternative scaling factor and we compared the two approaches for $Re \in  \{100,400,1000,2500,3200\}$. For this test, we cover a range of Reynolds numbers only in the laminar regime for the mere purpose of comparison with the $h_T$ scaling in terms of flow accuracy: indeed, the anisotropic scaling lacks convergence in the Newton solver for higher Reynolds numbers (see Table~\ref{tab:cavity:l2_errors:newARM}). 
We report in Table~\ref{tab:cavity:l2_errors:ARM:hF} and Table~\ref{tab:cavity:l2_errors:ARM:hT} the relative $\ell^2$-errors for low Reynolds numbers with anisotropic meshes using the scaling $h^\star_{\face}$ and $\hCell$, respectively. We compare with Ghia et al.\ \cite{GhiaGhiaShin1982} and Erturk et al.\ \cite{Erturk2009} depending on the available data (we explicitly indicate them in the table).
The use of $\hCell$ shows better behavior to reach the Newton convergence, caused by the overestimation of the Eddy viscosity.
However, smaller errors are obtained using  $h^\star_{\face}$ scaling in comparison with the results using the $\hCell$ scaling, due the overdiffusion of the latter model. 
We do not exclude that a more tailored guess choice would help the Newton convergence for the $h^\star_{\face}$ scaling; however, this investigation went beyond the goal of the contribution.
\begin{table}[H]
    \centering
    \caption{\emph{Lid-driven cavity flow}. Relative $\ell^2$-errors \reviewerAddA{for the velocity variable} with various $Re$ and ARM with scaling $h^\star_{\face}$. We compare with Ghia et al. \cite{GhiaGhiaShin1982} and Erturk et al. \cite{Erturk2009} depending on the available data.}
    \resizebox{\textwidth}{!}{
    \begin{tabular}{crrrrrr}
    \hline
    &\multicolumn{6}{c}{$Re$}\\ \cline{2-7}
    Scale& \multicolumn{1}{c}{100 \cite{GhiaGhiaShin1982}}  & \multicolumn{1}{c}{400 \cite{GhiaGhiaShin1982}}  & \multicolumn{1}{c}{1000 \cite{GhiaGhiaShin1982}}&  \multicolumn{1}{c}{1000 \cite{Erturk2009}}  & \multicolumn{1}{c}{2500 \cite{Erturk2009}}  & \multicolumn{1}{c}{3200 \cite{GhiaGhiaShin1982}}\\ \hline 
    Coarse &  1.016E-02&  8.209E-02&   1.659E-02&    2.720E-02&  - &-\\
    Medium &  1.078E-02&  8.275E-02&   1.324E-02&   2.049E-02 &  - &- \\     
    Fine   &  1.150E-02& 8.4074E-02 &1.330E-02&   1.267E-02&  1.970E-02&  3.562E-02\\ 
   \hline
    \end{tabular}}
    \label{tab:cavity:l2_errors:ARM:hF}
\end{table}

\begin{table}[H]
    \centering
    \caption{\emph{Lid-driven cavity flow}. Relative $\ell^2$-errors \reviewerAddA{for the velocity variable} for various $Re$ and ARM with scaling $h_T$. We compare with Ghia et al. \cite{GhiaGhiaShin1982} and Erturk et al. \cite{Erturk2009} depending on the available data. }
    \resizebox{\textwidth}{!}{
    \begin{tabular}{crrrrrr}
    \hline
    &\multicolumn{6}{c}{$Re$}\\ \cline{2-7}
    Scale& \multicolumn{1}{c}{100\cite{GhiaGhiaShin1982}} & \multicolumn{1}{c}{400\cite{GhiaGhiaShin1982}} & \multicolumn{1}{c}{1000\cite{GhiaGhiaShin1982}}&  \multicolumn{1}{c}{1000\cite{Erturk2009}}  & \multicolumn{1}{c}{2500\cite{Erturk2009}}  & \multicolumn{1}{c}{3200\cite{GhiaGhiaShin1982}}\\ \hline 
 Coarse&  1.449E-02&   8.577E-02&   1.178E-01&   9.429E-02&  - & -\\
 Medium&  1.363E-02&   8.413E-02&   8.823E-02&   6.972E-02&   1.602E-02&  - \\      
 Fine	     &  1.288E-02 &   5.055E-02   & 3.906E-02 &  1.054E-02 & 1.054E-02 & 1.118E-02 \\                   
   \hline
    \end{tabular}
    }
    \label{tab:cavity:l2_errors:ARM:hT}
\end{table} 
\end{remark}
\section{Conclusions}
\label{sec:conc}
We proposed the combination of the VEM with the discrete Smagorinsky model to solve convection-dominated NSE. We presented successful numerical tests that show the method's convergence to analytical solutions as the meshsize tends to zero. Moreover, we performed lid-driven cavity tests with different mesh families to assess the potential of VEM approximation for the Smagorinsky model. 
The IMH mesh, which is easily handled by the VEM, proved to be the best choice in terms of accuracy, efficiency, and robustness of the nonlinear solver when compared to uniform and anisotropic meshes. We also apply a different anisotropic scaling for the Smagorinsky term, i.e., based on the minimum face diameter. This variant was shown to be more accurate and to mitigate overdiffusion. Nonetheless, the scarce convergence behavior of the Newton solver with anisotropic meshes did not improve. \\
We remark that, to the best of our knowledge, this is the first work dealing with VEM approximation in turbulence modelling. For these reasons, many are the unexplored yet interesting paths of this research. For example, we aim at applying more problem-guided mesh refinement, based on \emph{ad hoc} a posteriori error estimator and/or based on the aspect ratio of the polygons \reviewerB{, and at providing tailored a priori error estimations for the VEM Smagorinsky model, with no restriction on the meshsize}. 
\section*{Acknowledgments}
\emergencystretch=3em 
This manuscript reflects only the authors’ views and opinions, and the Ministry cannot be considered responsible for them.
The author S.B.\ kindly acknowledges the funding by the European Union through project Next Generation EU, M4C2, PRIN 2022 PNRR project P2022BH5CB\_001 ``Polyhedral Galerkin methods for engineering applications to improve disaster risk forecast and management: stabilization-free operator-preserving methods and optimal stabilization methods''. The author K.L.C. acknowledges financial support by INdAM-GNCS
through the projects 2024 (CUP: E53C23001670001) and 2025 (CUP: E53C24001950001), and the PNRR M4C2 project of CN00000013 National Centre for HPC, Big Data and Quantum Computing (HPC) (CUP: E13C22000990001). Authors E.D. and S.R. acknowledge the Spanish Government - FEDER Fund Project PID2021-123153OBC21 funded by
MCIN/AEI/10.13039/501100011033/FEDER, UE. 
M.S. acknowledges the ``20227K44ME - Full and Reduced order modelling of coupled systems: focus on non-matching methods and automatic learning (FaReX)" project – funded by European Union – Next Generation EU  within the PRIN 2022 program (D.D. 104 - 02/02/2022 Ministero dell’Università e della Ricerca). Moreover, M.S. thanks the INdAM - GNCS Project ``Metodi numerici efficienti per problemi accoppiati in sistemi
complessi” (CUP\_E53C24001950001) and the ECCOMAS EYIC Grant ``CRAFT: Control and Reg-reduction in Applications for Flow Turbulence".
The author F.V. acknowledges the financial support by INdAM-research group GNCS, project title: Metodi numerici avanzati per equazioni alle derivate parziali in fenomeni di trasporto e diffusione - CUP E53C24001950001.


\bibliographystyle{abbrv}
\bibliography{bib.bib}

@article{GhiaGhiaShin1982,
    title = {{High-Re solutions for incompressible flow using the Navier-Stokes equations and a multigrid method}},
    journal = {Journal of Computational Physics},
    volume = {48},
    number = {3},
    pages = {387-411},
    year = {1982},
    issn = {0021-9991},
    doi = {https://doi.org/10.1016/0021-9991(82)90058-4},
    url = {https://www.sciencedirect.com/science/article/pii/0021999182900584},
    author = {U Ghia and K.N Ghia and C.T Shin},
}

@article{Erturk2009,
    author = {Erturk, Ercan},
    title  = {Discussions on driven cavity flow},
    journal = {International Journal for Numerical Methods in Fluids},
    volume = {60},
    number = {3},
    pages  = {275-294},
    doi = {https://doi.org/10.1002/fld.1887},
    url = {https://onlinelibrary.wiley.com/doi/abs/10.1002/fld.1887},
    eprint = {https://onlinelibrary.wiley.com/doi/pdf/10.1002/fld.1887},
    year = {2009}
}

@book{quarteroni2008numerical,
    title   =   {Numerical approximation of partial differential equations},
    author  =   {Quarteroni, Alfio and Valli, Alberto},
    volume  =   {23},
    year    =   {2008},
    publisher   =   {Springer Science \& Business Media}
}

@article{EnriqueROM,
    author = {Rebollo, Tom\'{a}s Chac\'{o}n and \'{A}vila, Enrique Delgado and M\'{a}rmol, Macarena G\'{o}mez and Ballarin, Francesco and Rozza, Gianluigi},
    title = {{On a Certified Smagorinsky Reduced Basis Turbulence Model}},
    journal = {SIAM Journal on Numerical Analysis},
    volume = {55},
    number = {6},
    pages = {3047-3067},
    year = {2017},
    doi = {10.1137/17M1118233}
}

@article{SamueleThesis,
  title={Numerical modeling of turbulence by Richardson numberbased and {VMS} models},
  author={Rubino, Samuele},
    journal={Ph.D. thesis},
  year={2014}
}

@article{Haferssas2018,
  author       = {Haferssas, Ryadh and Jolivet, Pierre and Rubino, Samuele},
  title        = {{Efficient and scalable discretization of the Navier–Stokes equations with LPS modeling}},
  journal      = {Computer Methods in Applied Mechanics and Engineering},
  volume       = {333},
  pages        = {371--394},
  year         = {2018},
  doi          = {10.1016/j.cma.2018.01.026},
}

@article{BeiraoBrezziCangianiEtAl2013,
    author = {Beir\~{a}o da Veiga, L. and Brezzi, F. and Cangiani, A. and Manzini, G. and Marini, L. D. and Russo, A.},
    title = {BASIC PRINCIPLES OF VIRTUAL ELEMENT METHODS},
    journal = {Mathematical Models and Methods in Applied Sciences},
    volume = {23},
    number = {01},
    pages = {199-214},
    year = {2013},
    doi = {10.1142/S0218202512500492},
    URL = {https://doi.org/10.1142/S0218202512500492},
}

@article{BeiraoBrezziMariniEtAl2014,
    author = {Beir\~{a}o da Veiga, L. and Brezzi, F. and Marini, L. D. and Russo, A.},
    title = {The Hitchhiker's Guide to the Virtual Element Method},
    journal = {Mathematical Models and Methods in Applied Sciences},
    volume = {24},
    number = {08},
    pages = {1541-1573},
    year = {2014},
    doi = {10.1142/S021820251440003X},
    URL = {https://doi.org/10.1142/S021820251440003X},
}

@article{VEM_Review2023,
    title = {The virtual element method}, 
    volume = {32}, 
    DOI = {10.1017/S0962492922000095}, 
    journal = {Acta Numerica}, 
    author = {Beirão Da Veiga, Lourenço and Brezzi, Franco and Marini, L. Donatella and Russo, Alessandro}, 
    year = {2023}, 
    pages = {123–202}
}

@article{BeiraoLovadinaVacca2017,
     author = {Beir\~ao da Veiga, L. and Lovadina, C. and Vacca, G.},
     title = {{Divergence free virtual elements for the Stokes problem on polygonal meshes}},
     journal = {ESAIM: Mathematical Modelling and Numerical Analysis },
     pages = {509--535},
     publisher = {EDP-Sciences},
     volume = {51},
     number = {2},
     year = {2017},
     doi = {10.1051/m2an/2016032},
     mrnumber = {3626409},
     zbl = {1398.76094},
     language = {en},
     url = {https://www.numdam.org/articles/10.1051/m2an/2016032/}
}

@article{BeiraoLovadinaVacca2018,
    author = {Beir\~ao da Veiga, L. and Lovadina, C. and Vacca, G.},
    title = {{Virtual Elements for the Navier--Stokes Problem on Polygonal Meshes}},
    journal = {SIAM Journal on Numerical Analysis},
    volume = {56},
    number = {3},
    pages = {1210-1242},
    year = {2018},
    doi = {10.1137/17M1132811},
    
    URL = {https://doi.org/10.1137/17M1132811
    },
    eprint = {https://doi.org/10.1137/17M1132811}
    }

@article{Brenner2018,
  author =       {Brenner, S. C. and Sung, L.-Y.},
  title =        {Virtual element methods on meshes with small edges or faces},
  journal =      {Mathematical Models and Methods in Applied Sciences},
  volume =       28,
  number =       07,
  pages =        {1291-1336},
  year =         2018,
  doi =          {10.1142/S0218202518500355},
}

@article{DassiVacca2020,
    title = {Bricks for the mixed high-order virtual element method: Projectors and differential operators},
    journal = {Applied Numerical Mathematics},
    volume = {155},
    pages = {140-159},
    year = {2020},
    note = {Structural Dynamical Systems: Computational Aspects held in Monopoli (Italy) on June 12-15, 2018.},
    issn = {0168-9274},
    doi = {https://doi.org/10.1016/j.apnum.2019.03.014},
    url = {https://www.sciencedirect.com/science/article/pii/S0168927419300674},
    author = {F. Dassi and G. Vacca},
}

@article{Cicuttin2025,
    title = {An implementation detail about the scaling of monomial bases in polytopal finite element methods},
    journal = {Applied Mathematics Letters},
    volume = {159},
    pages = {109281},
    year = {2025},
    issn = {0893-9659},
    doi = {https://doi.org/10.1016/j.aml.2024.109281},
    url = {https://www.sciencedirect.com/science/article/pii/S089396592400301X},
    author = {M. Cicuttin},
}

@article{Smagorinsky1963,
  title={General circulation experiments with the primitive equations: I. The basic experiment},
  author={Smagorinsky, Joseph},
  journal={Monthly Weather Review},
  volume={91},
  number={3},
  pages={99--164},
  year={1963},
  publisher={American Meteorological Society},
  doi={10.1175/1520-0493(1963)091<0099:GCEWTP>2.3.CO;2}
}

@book{chacon2014mathematical,
  title     = {Mathematical and Numerical Foundations of Turbulence Models and Applications},
  author    = {Chacon Rebollo, Tomas and Lewandowski, Roger},
  year      = {2014},
  publisher = {Springer New York},
  address   = {New York, NY},
  series    = {Modeling and Simulation in Science, Engineering and Technology},
  doi       = {10.1007/978-1-4939-0050-6},
  isbn      = {978-1-4939-0049-0}
}

@book{sagaut2001large,
  title     = {Large Eddy Simulation for Incompressible Flows: An Introduction},
  author    = {Sagaut, Pierre},
  year      = {2001},
  publisher = {Springer-Verlag Berlin Heidelberg},
  address   = {Berlin, Heidelberg},
  edition   = {1st},
  isbn      = {978-3-540-42204-6},
  doi       = {10.1007/978-3-662-04283-3}}

@misc{borio2024generalordervirtualelement,
      title={{General Order Virtual Element Methods for Neumann Boundary Optimal Control Problems in Saddle Point Formulation}}, 
      author={Andrea Borio and Francesca Marcon and Maria Strazzullo},
      year={2024},
      eprint={2411.08497},
      archivePrefix={arXiv},
      primaryClass={math.NA},
      url={https://arxiv.org/abs/2411.08497}, 
}

@article{Beirao2015a,
  title =        {A Virtual Element Method for elastic and inelastic problems on polytope meshes},
  journal =      {Computer Methods in Applied Mechanics and Engineering},
  volume =       295,
  pages =        {327-346},
  year =         2015,
  issn =         {0045-7825},
  doi =          {10.1016/j.cma.2015.07.013},
  author =       {{Beir\~{a}o da Veiga}, L. and Lovadina, C. and Mora, D.},
}

@article{Benedetto2016c,
  author =       {M. F. Benedetto and S. Berrone and A. Borio and S. Pieraccini and S. Scial\`{o}},
  title =        {A Hybrid Mortar Virtual Element Method For Discrete Fracture Network Simulations},
  journal =      {Journal of Computational Physics},
  fjournal =     {J. Comput. Phys.},
  year =         2016,
  volume =       306,
  pages =        {148--166},
  DOI =          {10.1016/j.jcp.2015.11.034},
}

@article{Borio2021,
  author =       {Borio, A. and Hamon, F. P. and Castelletto, N. and White, J. A. and Settgast,
                  R. R.},
  title =        {Hybrid mimetic finite-difference and virtual element formulation for coupled
                  poromechanics},
  journal =      {Computer Methods in Applied Mechanics and Engineering},
  volume =       383,
  pages =        113917,
  year =         2021,
  issn =         {0045-7825},
  doi =          {10.1016/j.cma.2021.113917}
}

@article{Dassi2021b,
  author =       {Dassi, F. and Lovadina, C. and Visinoni, M.},
  title =        {Hybridization of the virtual element method for linear elasticity problems},
  journal =      {Mathematical Models and Methods in Applied Sciences},
  volume =       31,
  number =       14,
  pages =        {2979-3008},
  year =         2021,
  doi =          {10.1142/S0218202521500676}
}

@ARTICLE{Dassi2024150,
  author =       {Dassi, F. and Mora, D. and Reales, C. and Velásquez, I.},
  title =        {{Mixed variational formulations of virtual elements for the polyharmonic operator
                  $(-\Delta)^n$}},
  journal =      {Computers and Mathematics with Applications},
  year =         {2024},
  volume =       {158},
  pages =        {150-166},
  doi =          {10.1016/j.camwa.2024.01.013}
}

@ARTICLE{Kumar2024,
  author =       {Kumar, S. and Mora, D. and Ruiz-Baier, R. and Verma, N.},
  title =        {Numerical Solution of the Biot/Elasticity Interface Problem Using Virtual Element
                  Methods},
  journal =      {Journal of Scientific Computing},
  year =         {2024},
  volume =       {98},
  number =       {3},
  doi =          {10.1007/s10915-023-02444-7},
  art_number =   {53}
}

@ARTICLE{Xu2024,
	author = {Xu, Bing-Bing and Fan, Wei-Long and Wriggers, Peter},
	title = {High-order {3D} virtual element method for linear and nonlinear elasticity},
	year = {2024},
	journal = {Computer Methods in Applied Mechanics and Engineering},
	volume = {431},
	doi = {10.1016/j.cma.2024.117258}
}

@book{roos2008robust,
  title={Robust Numerical Methods for Singularly Perturbed Differential Equations: 
         Convection-Diffusion-Reaction and Flow Problems.},
  author={Roos, H. G. and Stynes, M. and Tobiska, L.},
  volume={24},
  year={2008},
  series={Springer Series in Computational Mathematics},
  edition={second},
  publisher={Springer}
}

@article{MANZINI2022176,
title = {{Conforming virtual element approximations of the two-dimensional Stokes problem}},
journal = {Applied Numerical Mathematics},
volume = {181},
pages = {176-203},
year = {2022},
issn = {0168-9274},
doi = {https://doi.org/10.1016/j.apnum.2022.06.002},
url = {https://www.sciencedirect.com/science/article/pii/S0168927422001520},
author = {Gianmarco Manzini and Annamaria Mazzia}
}

@article{lilly1966representation,
  title={The representation of small-scale turbulence in numerical simulation experiments},
  author={Lilly, D},
  journal={NCAR Report},
  year={1966}
}

@article{kolmogorov1941dissipation,
  title={Dissipation of energy in the locally isotropic turbulence},
  author={Kolmogorov, Andrei Nikolaevich},
  journal={Proceedings of the Royal Society of London. Series A: Mathematical and Physical Sciences},
  volume={434},
  number={1890},
  pages={15--17},
  year={1991},
  publisher={The Royal Society London}
}

@article{kolmogorov1941local,
  title     =   {The local structure of turbulence in incompressible viscous fluid for very large {R}eynolds numbers},
  author    =   {Kolmogorov, Andrey Nikolaevich},
  journal   =   {Comptes Rendus de l'Acad\'emie des Sciences de l'URSS},
  volume    =   {30},
  pages     =   {301--305},
  year      =   {1941}
}

@book{fenics,
        title  = {Automated Solution of Differential Equations by the Finite Element Method},
        author = {Logg, A. and Mardal, K.A. and Wells, G.},
        year = {Berlin, 2012},
        publisher = {Springer-Verlag}
}

@ARTICLE{Strazzullo2025,
	author = {Strazzullo, Maria and Ballarin, Francesco and Iliescu, Traian and Rebollo, Tomás Chacón},
	title = {Variational multiscale evolve and filter strategies for convection-dominated flows},
	year = {2025},
	journal = {Computer Methods in Applied Mechanics and Engineering},
	volume = {438},
	doi = {10.1016/j.cma.2025.117811}
}

@ARTICLE{Zoccolan2025237,
	author = {Zoccolan, Fabio and Strazzullo, Maria and Rozza, Gianluigi},
	title = {{A Streamline Upwind Petrov-Galerkin Reduced Order Method for Advection-Dominated Partial Differential Equations Under Optimal Control}},
	year = {2025},
	journal = {Computational Methods in Applied Mathematics},
	volume = {25},
	number = {1},
	pages = {237 – 260},
	doi = {10.1515/cmam-2023-0171}
}

@ARTICLE{Zoccolan2024,
	author = {Zoccolan, Fabio and Strazzullo, Maria and Rozza, Gianluigi},
	title = {Stabilized weighted reduced order methods for parametrized advection-dominated optimal control problems governed by partial differential equations with random inputs},
	year = {2024},
	journal = {Journal of Numerical Mathematics},
	doi = {10.1515/jnma-2023-0006}
}

@article{Ivagnes2025,
  author       = {Ivagnes, A. and Strazzullo, M. and Girfoglio, M. and Iliescu, T. and Rozza, G.},
  title        = {Data-Driven Optimization for the Evolve-Filter-Relax Regularization of Convection-Dominated Flows},
  journal      = {International Journal for Numerical Methods in Engineering},
  volume       = {126},
  number       = {9},
  pages        = {e70042},
  year         = {2025},
  doi          = {10.1002/nme.70042},
  url          = {https://doi.org/10.1002/nme.70042}
}

@article{Khani_Waite_2015, 
    title={Large eddy simulations of stratified turbulence: the dynamic Smagorinsky model}, 
    volume={773}, 
    DOI={10.1017/jfm.2015.249}, 
    journal={Journal of Fluid Mechanics}, 
    author={Khani, Sina and Waite, Michael L.}, 
    year={2015}, 
    pages={327–344}
}

@article{Rozema2022,
  title = {Local dynamic gradient Smagorinsky model for large-eddy simulation},
  author = {Rozema, Wybe and Bae, H. Jane and Verstappen, Roel W. C. P.},
  journal = {Phys. Rev. Fluids},
  volume = {7},
  issue = {7},
  pages = {074604},
  numpages = {26},
  year = {2022},
  month = {Jul},
  publisher = {American Physical Society},
  doi = {10.1103/PhysRevFluids.7.074604},
}

@article{Winckelmans2001,
    author = {Winckelmans, Grégoire S. and Wray, Alan A. and Vasilyev, Oleg V. and Jeanmart, Hervé},
    title = {Explicit-filtering large-eddy simulation using the tensor-diffusivity model supplemented by a dynamic Smagorinsky term},
    journal = {Physics of Fluids},
    volume = {13},
    number = {5},
    pages = {1385-1403},
    year = {2001},   
   doi = {10.1063/1.1360192}
}

@book{berselli2006mathematics,
  title={Mathematics of large eddy simulation of turbulent flows},
  author={Berselli, Luigi C and Iliescu, Traian and Layton, William J},
  year={2006},
  publisher={Springer}
}

@article{rebollo2014numerical,
  title={Numerical approximation of the Smagorinsky turbulence model applied to the primitive equations of the ocean},
  author={Rebollo, Tom{\'a}s Chac{\'o}n and Hecht, Fr{\'e}d{\'e}ric and M{\'a}rmol, Macarena G{\'o}mez and Orzetti, Giordano and Rubino, Samuele},
  journal={Mathematics and Computers in Simulation},
  volume={99},
  pages={54--70},
  year={2014},
  publisher={Elsevier}
}

@article{germano1991dynamic,
  title={A dynamic subgrid-scale eddy viscosity model},
  author={Germano, Massimo and Piomelli, Ugo and Moin, Parviz and Cabot, William H},
  journal={Physics of fluids a: Fluid dynamics},
  volume={3},
  number={7},
  pages={1760--1765},
  year={1991},
  publisher={American Institute of Physics}
}

@article{iliescu2003numerical,
  title={A numerical study of a class of LES models},
  author={Iliescu, T and John, V and Layton, WJ and Matthies, G and Tobiska, L},
  journal={International Journal of Computational Fluid Dynamics},
  volume={17},
  number={1},
  pages={75--85},
  year={2003},
  publisher={Taylor \& Francis}
}

@article{lilly1992proposed,
  title={A proposed modification of the Germano subgrid-scale eddy viscosity model},
  author={Lilly, DK},
  journal={Phys. Fluids A},
  volume={4},
  number={3},
  pages={633--635},
  year={1992}
}

@article{pare1992existence,
  title={Existence, uniqueness and regularity of solution of the equations of a turbulence model for incompressible fluids},
  author={Par{\'e}; s, Carlos},
  journal={Applicable Analysis},
  volume={43},
  number={3-4},
  pages={245--296},
  year={1992},
  publisher={Taylor \& Francis}
}

@article{Sorgente2021,
  author    = {T. Sorgente and S. Biasotti and G. Manzini and M. Spagnuolo},
  title     = {The role of mesh quality and mesh quality indicators in the virtual element method},
  journal   = {Advances in Computational Mathematics},
  year      = {2021},
  volume    = {48},
  number    = {1},
  pages     = {3},
  doi       = {10.1007/s10444-021-09913-3},
  url       = {https://doi.org/10.1007/s10444-021-09913-3},
  issn      = {1572-9044}
}

@article{berrone2025polydim,
title = {{POLYDIM: A C++ library for POLYtopal DIscretization Methods}},
journal = {Computer Physics Communications},
volume = {320},
pages = {109937},
year = {2026},
issn = {0010-4655},
doi = {10.1016/j.cpc.2025.109937},
author = {Stefano Berrone and Andrea Borio and Gioana Teora and Fabio Vicini},
}

\end{document}